\documentclass[10pt]{article}
\usepackage{pb-diagram}
\usepackage{amsmath}
\usepackage{amssymb}
\usepackage{amscd}

\newtheorem{Df}{Definition}[section]
\newtheorem{Te}[Df]{Theorem}
\newtheorem{Po}[Df]{Proposition}
\newtheorem{Cr}[Df]{Corollary}
\newtheorem{Lm}[Df]{Lemma}
\newtheorem{Ca}[Df]{Claim}
\newtheorem{Cn}[Df]{Conjecture}

\newcommand{\Bdf}{\begin{Df}}
\newcommand{\Edf}{\end{Df}}
\newcommand{\Bte}{\begin{Te}}
\newcommand{\Ete}{\end{Te}}
\newcommand{\Bpo}{\begin{Po}}
\newcommand{\Epo}{\end{Po}}
\newcommand{\Bcr}{\begin{Cr}}
\newcommand{\Ecr}{\end{Cr}}
\newcommand{\Blm}{\begin{Lm}}
\newcommand{\Elm}{\end{Lm}}
\newcommand{\Bca}{\begin{Ca}}
\newcommand{\Eca}{\end{Ca}}
\newcommand{\Bcn}{\begin{Cn}}
\newcommand{\Ecn}{\end{Cn}}
\newcommand{\Bdm}{{\it Proof.}\ }
\newcommand{\Edm}{\rule{2mm}{2mm}}
\newcommand{\Rm}{{\it Remark \arabic{section}.\arabic{Df} \ }}
\newcommand{\Ea}{{\it Example \arabic{section}.\arabic{Df} \ }}

\begin{document}

\title{\bf{Calabi-Yau algebras viewed as deformations of Poisson algebras}}
\author{\bf{Roland Berger and Anne Pichereau}\footnote{The second author was supported by the french ANR projet
ANR-09-RPDOC-009-01.}\\ Universit\'e de Lyon, Institut Camille Jordan (UMR 5208) \\ Universit\'e de Saint-Etienne, Facult\'{e} des Sciences\\ 23, Rue Docteur Paul Michelon, 42023 Saint-Etienne
Cedex 2,
France \\
\emph{E-mail}: Roland.Berger@univ-st-etienne.fr,\\ Anne.Pichereau@univ-st-etienne.fr}
\date{}
\maketitle
\begin{abstract}
From any algebra $A$ defined by a single non-degenerate homogeneous quadratic
relation $f$, we prove that the quadratic algebra $B$ defined by the potential $w=fz$ is 3-Calabi-Yau. The algebra $B$ can be viewed as a 3-Calabi-Yau completion of Keller. The algebras
$A$ and $B$ are both Koszul. The classification of the algebras $B$ in three generators, i.e.,
when $A$ has two generators, leads to three types of algebras. The second type (the most
interesting one) is viewed as a deformation of a Poisson algebra $S$ whose Poisson bracket is
non-diagonalizable quadratic.
Although the potential of $S$ has non-isolated singularities, the homology of $S$ is computed. Next the Hochschild homology of $B$ is obtained.

\end{abstract}

\section{Introduction}
The deformation theory offers a way to study non-commutative algebras $A$ by
examining associated Poisson algebras $S$. The idea is to get (more or less explicitly) invariants
attached to $S$ by the use of differential calculus on the commutative algebra $S$, and then to
deduce from them some invariants attached to $A$. Invariants of interest in various domains are of
homological nature. In particular, a natural question is the following (there is the same in cohomology).

(Q) Is the Hochschild homology HH$_{\bullet}(A)$ of $A$ isomorphic to the Poisson
homology HP$_{\bullet}(S)$ of $S$?
 A positive answer was given by Kontsevich~\cite{kont:quant} in cohomology when $S$
is the algebra of the C$^{\infty}$ functions on $\mathbb{R}^n$ endowed with any Poisson bracket $\pi$ extended to the
space of formal series $S[[\hbar]]$, and $A$ is the space $S[[\hbar]]$ whose usual commutative product is replaced by the
Kontsevich star product $\star _{\pi}$ (see e.g.~\cite{cktb:def} for further developments). The question (Q) was initiated by Brylinski~\cite{bry:poi} when the algebra $A$ is
filtered such that gr$(A)$ is assumed to be commutative and smooth. Then gr$(A)$ is naturally a Poisson
algebra $S$, and there is the Brylinski spectral sequence
\begin{equation} \label{ssbry}
E^2= HP_{\bullet}(S) \Longrightarrow HH_{\bullet}(A).
\end{equation} In this context, we can replace the question (Q) by the following one.

(Q') Does the Brylinski spectral sequence degenerate at $E^2$?

As shown by Kassel~\cite{kas:homcy}, this is the case if $S$ is any polynomial algebra
whose generators have all degree 1, so that the Poisson bracket of $S$ is of degree $\leq 1$ and $A$
is a Sridharan enveloping algebra. Actually, Kassel proved more precisely that the symmetrization
defines an isomorphism from the Poisson complex of $S$ to the Koszul complex of $A$. The latter
result was generalized (in a weaker form) in~\cite{hot:orbi} to the crossed products of enveloping algebras.

Van den Bergh showed that the question (Q') has a positive answer if $A$ is a
generic Sklyanin algebra in three generators~\cite{vdb:Kth}. In this case, the Poisson bracket of $S$ (which
is quadratic) is derived from a ``Poisson potential'' $\phi$ such that

(IS) the origin is an isolated singularity of the polynomial $\phi$.

In the paper of Van den Bergh, the property (IS) appears as essential in the
computation of HP$_{\bullet}(S)$, and it is the same for the other examples (quadratic, cubic, or more generally quasi-homogeneous)
considered further~\cite{nic:cubic, pic:hom, pel:four}. In this paper, we present an example of quadratic algebra, called $B$ in the text, for which the previous programme is
performed without the property (IS). In other words, in our example of quadratic algebra $B$, the
Poisson bracket of $S$ (which is non-diagonalizable quadratic) is derived from a potential $\phi$ having
non-isolated singularities. Nevertheless, HP$_{\bullet}(S)$ is explicitly computed (Section 5 below). Next we
prove that any Poisson cycle can be lifted to a Koszul cycle (since $B$ is Koszul, it is more convenient
to use Koszul complex, instead of Hochschild complex, in order to compute HH$_{\bullet}(B)$). Finally we get the computation of HH$_{\bullet}(B)$ as stated in the following.
\Bte \label{calc}
Let $B$ be the $\mathbb{C}$-algebra defined by the generators $x$, $y$, $z$, and the
following relations $$zy=yz+2xz, \ zx=xz, \ yx=xy + x^2.$$ Let $S$ be the polynomial $\mathbb{C}$-algebra in $x$, $y$, $z$, endowed with the Poisson bracket derived from the potential $\phi = -x^2 z$.
Then the Hochschild homology of $B$ is isomorphic by a specific canonical morphism (an edge morphism in the Brylinski spectral sequence) to the Poisson homology of $S$ and
is given by \begin{eqnarray*}
HH_0(B) &\simeq& x\mathbb{C}[y]\oplus \mathbb{C}[y,z];\\
 {}\\
HH_1(B) &\simeq& \mathbb{C}[\phi]
\left(\begin{smallmatrix}xz\\0\\-x^2\end{smallmatrix}\right)\oplus\mathbb{C}[z]
\left(\begin{smallmatrix}z^2\\0\\-xz\end{smallmatrix}\right)\oplus  \renewcommand{\arraystretch}{0.7}
  \bigoplus_{\begin{array}{c}\scriptstyle n\in \mathbb N^*\\\scriptstyle 0\leq k\leq
n\end{array}}
 \mathbb{C}  \left(\begin{smallmatrix}0\\ky^{k-1}z^{n-k}\\(n-k)y^kz^{n-1-k}\end{smallmatrix}\right)\\
&\oplus &\bigoplus_{n\in \mathbb
N}\,\mathbb{C}\left(\begin{smallmatrix}y^n\\n\,xy^{n-1}\\0\end{smallmatrix}\right)
\oplus  \renewcommand{\arraystretch}{0.7}
  \bigoplus_{\begin{array}{c}\scriptstyle n\in \mathbb N\\ \scriptstyle 1\leq k\leq
n+1\end{array}}\mathbb C\left(\begin{smallmatrix}(2n+3)\, yz\\-3k\,xz\\
(-2n+3(k-1))\,xy\end{smallmatrix}\right)y^{k-1}z^{n+1-k};
\end{eqnarray*}
\begin{eqnarray*}
HH_2(B) &\simeq& \mathbb C[\phi]
\left(\begin{smallmatrix}x\\y\\z\end{smallmatrix}\right)\,\oplus\,   \left(x\mathbb C[\phi]\oplus z\mathbb
C[z]\right)\left(\begin{smallmatrix}0\\1\\0\end{smallmatrix}\right)\\
 &&\qquad \qquad \qquad \qquad \qquad \qquad \oplus \renewcommand{\arraystretch}{0.7}
  \bigoplus_{\begin{array}{c}\scriptstyle n\in \mathbb N\\\scriptstyle 0\leq k\leq
n\end{array}}\mathbb C\left(\begin{smallmatrix}(k+1) x\\ (2(n-k)+1)\,y\\
-2(k+1)\,z\end{smallmatrix}\right)y^{k}z^{n-k};\\
  {}\\
  HH_3(B) &\simeq& \mathbb C[\phi];\\
  HH_p(B) &\simeq& 0 \ \mathrm{if} \ p\geq 4.
\end{eqnarray*}
\Ete

Since $B$ is 3-Calabi-Yau, the Hochschild cohomology is then immediate:
$HH^{\bullet}(B)=HH_{3-\bullet}(B)$. Remark that the duality $HP^{\bullet}(S)=HP_{3-\bullet}(S)$ holds since the Poisson bracket derives from a potential.

Our example of Koszul algebra $B$ belongs to a large class of Koszul 3-Calabi-Yau
algebras still denoted by $B$. These algebras $B$ can be viewed as 3-Calabi-Yau completions of algebras $A$ of global dimension 2 (consequently they are quasi-isomorphic to Ginzburg algebras~\cite{vg:cy}). Calabi-Yau completions were defined and studied by Keller in~\cite{bk:compl} (with an appendix by Van den Bergh). In particular we refer to Section 6.8 of~\cite{bk:compl}. The Calabi-Yau property of our algebras $B$ is then immediate from Theorem 4.8 in~\cite{bk:compl}. Let us explain how we define these algebras $B$ (see Section 2 below for details). For any
non-degenerate quadratic relation $f=\sum_{i,j} f_{ij} x_i x_j$ in $n \geq 2$ non-commutative
variables $x_i$, the algebra $A$ defined by the single relation $f$ is Koszul and AS-Gorenstein of
global dimension 2~\cite{dv:multi}, and $A$ is Calabi-Yau if and only if $f$ is symplectic~\cite{rb:gera}.
Add an extra generator $z$ to the $x_i$'s and consider the quadratic algebra $B$
defined by the potential $w=fz$ (see e.g.~\cite{vg:cy, bock:graded, rbrt:pbw} for algebras defined by
potentials). The properties of $B$ that we shall obtain in Section 2 work out over any field and are stated in the
following.
\Bte \label{clas}
Let $f=\sum_{i,j} f_{ij} x_i x_j$ be any non-degenerate quadratic relation over a
field $k$ in non-commutative variables $x_1, \ldots , x_n$ ($n \geq 2$) and let $z$ be an extra
generator. Let $A$ be the $k$-algebra defined by the generators $x_1, \ldots , x_n$ and the
single relation $f$. Let $B$ be the $k$-algebra defined by the generators $x_1, \ldots , x_n , z$ and the
potential $w=fz$. The following hold \\ 1) $B$ is a skew polynomial ring over $A$ in the generator $z$ and defined by an
automorphism of $A$.
\\ 2) $z$ is normal in $B$, that is $Bz=zB$.
\\ 3) If $f$ is alternating, then $z$ is central in $B$. The converse holds if the
characteristic of $k$ is $\neq 2$.
\\ 4) $B$ is Koszul and 3-Calabi-Yau. Moreover, $B$ is a domain.
\\ 5) The Hilbert series $h_B(t)$ of the graded algebra $B$ is $$h_B(t)=(1-(n+1)t+(n+1)t^2 - t^3)^{-1}.$$
5) The Gelfand-Kirillov dimension GK.dim$(B)$ of $B$ is finite if and only if $n=2$, and in this case GK.dim$(B)=3$.
\\ 6) $B$ is left (or right) noetherian if and only if $n=2$.
\Ete

When $n=2$, the graded $\mathbb{C}$-algebras $B$ are classified in three types
(Section 4 below): the polynomial algebra in $x$, $y$ and $z$ (classical type), the algebra of Theorem \ref{calc} (Jordan type), a family of
quantum spaces in $x$, $y$, $z$ (quantum type). The second type is the one of interest for us, since it is
well-known that the question (Q) has a positive answer in the first or third type (if the quantum parameter $q$
is not a root of unity). Hochschild homology in the quantum type can be deduced from
Wambst's result (\cite{wam:kosz}, Th\'eor\`eme 6.1). When $n=2$, the $\mathbb{C}$-algebras $B$ are AS-regular of global dimension 3, and it is important to notice that their invariants $j$ (in terminology
of~\cite{as:regular}) are \emph{infinite}, unlike Sklyanin algebras. The first type and the third type can
be considered as limits of Sklyanin algebras by vanishing the parameter $c$ used in~\cite{vdb:Kth}, but such a process is not possible for the second type which stands really apart.
Recently, Smith has given a detailed study of a remarkable algebra having seven
generators and defined in terms of octonions~\cite{smith:octo}. To avoid confusion with our notation, let
us call $C$ the Smith algebra. Actually, $C$ does not belong to the class of algebras $B$ defined in Theorem \ref{clas},
since in characteristic zero $C$ has no normal element except the elements of $k$ (\cite{smith:octo}, Proposition 11.2).
However $C$ has many properties in common with our algebras $B$: $C$ is Koszul and 3-Calabi-Yau, $C$ is defined by an
explicit potential, $C$ is a skew polynomial ring over an algebra $A$ in the last generator (but $C$ is
defined by a derivation of $A$, not an automorphism), and the single relation of $A$ is symplectic in the first
six generators. Su\'{a}rez-Alvarez has obtained similar properties for a class of
algebras containing $C$ and defined from any oriented Steiner triple
system~\cite{msa:steiner}. It would be satisfactory to enlarge naturally the class
of algebras $B$ in order to include the
Su\'{a}rez-Alvarez algebras, and one would expect that for any AS-regular algebra there is a skew polynomial ring in one variable over it which is Calabi-Yau.
\\ \emph{Acknowledgements.} We would like to thank the referee for valuable suggestions and comments.

\setcounter{equation}{0}

\section{A family of 3-Calabi-Yau algebras}

A down-to-earth approach of non-commutative projective algebraic geometry consists
in studying graded algebras defined by generators (assumed to be of degree 1) satisfying some homogeneous
non-commutative polynomial relations. Following this na\"{\i}ve approach, the first class to study is certainly the class of
non-commutative quadrics, i.e., the class of non-commutative graded algebras defined by a single quadratic relation. It is a bit surprising that this class can be used as a toy model (see~\cite{rb:gera}) in order to introduce to
several duality theories (Koszul duality, AS-Gorenstein duality and Calabi-Yau duality) playing a basic role
in more sophisticated approaches.
Throughout the paper, the algebras of this class will be denoted by $A$, sometimes
$A(f)$ or $A(M)$ if we want to specify the na\"{\i}ve non-commutative quadric which is just a quadratic relation $f$ of matrix
$M$ in the non-commutative generators. Our first goal will be to show that such an algebra $A$ is the quotient of a
quadratic graded algebra $B$ defined by a potential depending on $f$. The study of the algebras $B=B(f)=B(M)$ is
the main purpose of this paper.
 Let us fix the notation. For $n\geq 1$, $k\langle x_1, \ldots , x_n\rangle$ denotes
the graded free associative algebra over a field $k$, generated by $x_1, \ldots , x_n$ assumed to be of degree 1. Let us give a non-zero element
$f=\sum _{1\leq i,j \leq n} f_{ij} x_i x_j$ homogeneous of degree 2 in this algebra, or equivalently a non-zero $n\times n$
matrix $M=(f_{ij})_{1\leq i,j \leq n}$ with entries in $k$. Then $A=A(f)=A(M)$ denotes the quadratic graded algebra defined as the quotient of
$k\langle x_1, \ldots , x_n\rangle$ by the two-sided ideal generated by $f$.
Let us recall the properties of $A$ (see~\cite{rb:gera, jz:nnoeth} for the proofs and for the
definitions of Koszul, AS-Gorenstein or Calabi-Yau algebras).
\Bpo \label{recall}
Let $A=A(f)=A(M)$ be as above. \\ 1) The quadratic algebra $A$ is Koszul. \\ 2) The global dimension of $A$ is equal to 2, except if $f$ is symmetric of rank
1 (in this case, the global dimension is infinite).
\\ 3) $A$ is AS-Gorenstein if and only if $f$ is non-degenerate (i.e., $M$ is invertible).
\\ 4) $A$ is 2-Calabi-Yau if and only if $f$ is non-degenerate and skew-symmetric.
\\ 5) If the global dimension of $A$ is equal to 2, the Hilbert series of the graded
algebra $A$ is given by $$h_A(t)= (1-nt +t^2)^{-1}.$$
Otherwise, one has $$h_A(t)= (1-nt + t^2 -t^3 + t^4 - \cdots )^{-1}.$$
6) The Gelfand-Kirillov dimension GK.dim$(A)$ of $A$ is equal to 0 if $n=1$, to
$\infty$ if $n >2$, and if $n=2$, it is equal to the global dimension.
\\ 7) If $f$ is non-degenerate, $A$ is left (or right) noetherian if and only if $n=1$ or $n=2$.
\\ 8) If $f$ is non-degenerate, $A$ is a domain if and only if $n\geq 2$.\Epo

Set $F=k\langle x_1, \ldots , x_n, z\rangle$ where $z$ is an extra generator of
degree 1. We refer to~\cite{vg:cy, bock:graded, rbrt:pbw} for more details on the definitions and basic properties concerning on algebras defined
by a potential.  The elements of $F_{cyc} = F/[F,F]$ are called \emph{potentials}. The $k$-vector space $F_{cyc}$ is sometimes identified to the space of the cyclic
sums $c(a)$ when $a$ runs over $F$. Let us define our potential as being the class $\overline{w}$ in $F_{cyc}$ of $w\in F$ where $$w=fz = \sum _{1\leq i,j \leq n} f_{ij}\, x_i x_j z,$$
or as being the cyclic sum $$c(w)=\sum _{1\leq i,j \leq n} f_{ij} (x_i x_j z + z x_i
x_j + x_j z x_i).$$
Let us denote by $B=B(f)=B(M)$ the algebra defined by this potential. This is the
quotient of the free algebra $F$ by the cyclic partial derivatives $\partial_{x_1}(w), \ldots , \partial_{x_n}(w),
\partial_z(w)$, where
\begin{equation}
\partial_{x_i}(w)= \sum_{1\leq j \leq n} (f_{ij}\,x_j z + f_{ji}\, z x_j),\ \ 1 \leq
i \leq n,
\end{equation}
\begin{equation}
\partial_{z}(w)= f.
\end{equation}

So $B$ is a quadratic graded algebra. The obvious morphism of graded algebras from $A$ to the
quotient of $B$ by $z$ is clearly an isomorphism. In all the following, we abbreviate this isomorphism by $A \cong B/(z)$. The following lemma
allows us to consider $A$ as a quadratic subalgebra of $B$.
\Blm \label{subalg}
The morphism of graded algebras $A \rightarrow B$ induced by the composite
$$k\langle x_1, \ldots , x_n\rangle \stackrel{can}{\longrightarrow} F
\stackrel{\mathrm{can}}{\longrightarrow} B$$
is injective.
\Elm
\Bdm
It suffices to observe that the ideal which defines the algebra $B$ is homogeneous with respect to the grading by $\deg_z$.
\Edm
\\

 Denote by $V_A$ and $R_A$ (respectively $V_B$ and $R_B$) the space of generators and
the space of relations of $A$ (resp. $B$). One has $V_A=kx_1 \oplus \ldots \oplus kx_n$, $R_A=kf$, $V_B=V_A \oplus kz$ and $R_B=
(\sum_{1\leq i \leq n} k\partial_{x_i}(w)) \oplus R_A$.
\Blm \label{basic} 1) We have $\dim (R_B) \geq \emph{rk}(M) + 1$.
\\ 2) In particular, if $f$ is non-degenerate, then $\dim
R_B =n+1$ and $\partial_{x_1}(w), \ldots , \partial_{x_n}(w)$ are $k$-linearly independent in $F$. \\ 3) If $f$ is non-degenerate, then $Az=zA$ and the element $z$ is normal in the algebra $B$, that is, $Bz=zB$.
\Elm
\Bdm
Write down the relations $\partial_{x_i}(w)= 0$, $\ 1 \leq i \leq n$, as the linear
system (in the space $B$)
\begin{equation} \label{syst1}
M \left( \begin{array}{c}
x_1 z \\ \vdots \\ x_n z \end{array} \right) = -\, ^t M \left( \begin{array}{c}
z x_1  \\ \vdots \\ z x_n \end{array} \right)
\end{equation}
with unknowns $x_1 z, \ldots , x_n z$ viewed in $B$. Reducing this system to a
triangular form provides $p=rk(M)$ equalities in $B$ beginning by $p$ distincts pivots $x_i z$, so that these equalities viewed as elements of $R_B$
are linearly independent. From this the inequality in 1) follows. This shows as well that if $M$ is invertible, the elements $x_1 z, \ldots , x_n z$
belong to $zA$, proving that $Az \subseteq zA$ and $Bz \subseteq zB$. The opposite inclusions are obtained similarly by reading the linear system (\ref{syst1}) from the right to
the left. \Edm
\\

It is also easy to deduce the inequality in 1) from the following relation
\begin{equation} \label{rang}
\dim(R_{B(M)})=\emph{rk} \left( \begin{array}{cc}M & \,^tM \end{array} \right) +1,
\end{equation}
which is immediate from the following
\begin{equation} \label{rel}
\left( \begin{array}{c}
\partial_{x_1}(w) \\ \vdots \\ \partial_{x_n}(w) \end{array} \right)
= \left( \begin{array}{cc}
M & \,^t M
\end{array} \right) \left( \begin{array}{c}
x_1 z \\ \vdots \\ z x_n \end{array} \right)
\end{equation}
where the latter column is formed by the $2n$ elements $x_1z, \ldots , x_n z, z x_1,
\ldots z x_n$. Note that as a consequence of (\ref{rang}), we have equality in 1) if $M$ is symmetric or antisymmetric (not
necessarily invertible), and that for $n=2$ and $f=x_1^2 + x_1 x_2$, the inequality is
strict.
\\ \\
\addtocounter{Df}{1}
\Ea \label{central}
Assume that $f$ is symplectic, i.e., non-degenerate and alternating. It follows that $n=2p$ is even
and we can choose generators $x_1, \ldots, x_p, y_1, \ldots , y_p$ such that $f=\sum_{1\leq i \leq p} (x_iy_i - y_i x_i)$. Then $\partial_{x_i}(w)=y_i z
- z y_i$ and $\partial_{y_i}(w)= z x_i - x_i z$ for $1\leq i \leq p$. In this case, $z$ is \emph{central} in $B$ and one has \begin{equation}
c(w)= \sum_{1\leq i \leq p} Ant(x_i, y_i, z)
\end{equation}
where $Ant(a,b,c)$ denotes the antisymmetrizer of $a$, $b$, $c$.
\\

We continue the study of the algebra $B=B(f)$. Throughout the rest of this
section, we assume that $f$ is \emph{non-degenerate}. Lemma \ref{basic} shows that $Az=zA$ is a sub-$A$-bimodule of $B$. Let
$T_A(Az)$ denote the tensor algebra of the $A$-bimodule $Az$ (on tensor algebras of bimodules, see
e.g.~\cite{ks:qg}, p.485). The $k$-algebra $T_A(Az)$ is generated by $A \oplus Az$, hence by $x_1, \ldots x_n, z$. Let $$u : T_A (Az) \rightarrow B$$
be the natural morphism of $k$-algebras extending the inclusion $A \oplus Az\hookrightarrow B$. The inclusions $x_1, \ldots , x_n \hookrightarrow A$ and $z\hookrightarrow Az$ define a morphism of
$k$-algebras $$F \rightarrow T_A (Az)$$
factoring out the relations of $B$ and providing the morphism of $k$-algebras $$v: B \rightarrow T_A (Az)$$
which leaves ``fixed'' $x_1, \ldots x_n, z$. Then $v\circ u$ and $u\circ v$ coincide with the respective identity maps since the $k$-algebras $T_A(Az)$ and
$B$ are both generated by $x_1, \ldots , x_n, z$. Thus $u$ and $v$ are isomorphisms, inverse to each other. The following lemma is
crucial.
\Blm \label{free1}
The left (resp. right) sub-$A$-module $Az$ (resp. $zA$) of $B$ is free generated by
$z$.
\Elm
\Bdm
Let $a \in k\langle x_1, \ldots , x_n\rangle$ be homogeneous of degree $\ell \geq 0$
such that $az$ belongs to the two-sided ideal of $F$ generated by $\partial_{x_1}(w), \ldots , \partial_{x_n}(w), f$. Our aim is to prove
that $a$ belongs to the two-sided ideal of $k\langle x_1, \ldots , x_n\rangle$ generated by $f$. Write
\begin{equation}
az= \sum_{1\leq i \leq n, \alpha ,\beta \in \cal{M}} \lambda_{i,\alpha, \beta}\
\alpha \partial_{x_i}(w) \beta \ + \ \sum_{\alpha, \beta \in \cal{M'} } \mu_{\alpha, \beta}\ \alpha f \beta,
\end{equation}
where $\cal{M}$ denotes the set of (non-commutative) monomials in $x_1, \ldots ,x_n$
and $\cal{M'}$ is the union of $\cal{M}$ with the singleton $\{ z \}$,
where $\lambda_{i,\alpha, \beta}$ and $\mu_{\alpha, \beta}$ are in $k$, and where
the sums are finitely supported. In the second sum, one has $\deg_z \alpha + \deg_z \beta =1$ where $\deg_z$ denotes the degree in $z$. Viewing
the linear system (\ref{syst1}) in the free algebra $F$ and reducing it to the triangular form, we see that $z$ in $\alpha$ or $\beta$ can be put on the right modulo the elements $\partial_{x_i}(w)$.
Moreover the non-commutative Euler relation \begin{equation} \label{euler}
\sum_{1\leq i \leq n} \partial_{x_i} (w) x_i + fz = \sum_{1\leq i \leq n} x_i \partial_{x_i} (w) + zf = c(w)
\end{equation}
associated to the potential $w$ shows that the $z$ appearing in each $\alpha f
\beta$ can be put completely on the right modulo the elements $\partial_{x_i}(w)$. So we can write
\begin{equation}
az= \sum_{1\leq i \leq n, \alpha ,\beta \in \cal{M}} \lambda_{i,\alpha, \beta}\
\alpha \partial_{x_i}(w) \beta \ + \ ( \sum_{\alpha, \beta \in \cal{M} } \mu_{\alpha, \beta}\ \alpha f \beta \,) z,
\end{equation}
by keeping the same notation for the coefficients. Denote by $S$ the first sum. In
order to conclude, it is sufficient to prove that $S=0$. We have
\begin{equation}
S= \sum_{1\leq i,j \leq n, \alpha ,\beta \in \cal{M}}( \lambda_{i,\alpha, \beta}\
f_{ij} \alpha x_j z \beta + \lambda_{i,\alpha, \beta}\ f_{ji} \alpha z x_j \beta),
\end{equation}
where $\deg \alpha + \deg \beta =\ell -1$ and $\deg$ denotes the total degree in
$x_1, \ldots , x_n$.

From $S=(a-\sum_{\alpha, \beta \in \cal{M} } \mu_{\alpha, \beta}\ \alpha f \beta \,)
z$, we are going to deduce that the coefficients $\lambda_{i,\alpha, \beta}$ vanish inductively. Choose $\beta $ with the maximal
degree $\ell-1$, so that $\alpha=1$. Since $\deg (x_j \beta) >\deg (\beta')$  for any $\alpha'$, $\beta'$ appearing in a term
$\alpha' x_k z \beta'$ such that $\deg \alpha' + \deg \beta'= \ell-1$, the coefficient of $zx_j \beta$ vanishes, hence $\sum_{1\leq i \leq n} \lambda_{i,1, \beta } f_{ji}=0$
for $j=1, \ldots , n$. As $M$ is invertible, this implies that $\lambda_{i,1, \beta }=0$ for $i=1, \ldots ,n$. Thus we can remove all the elements
$\lambda_{i,1, \beta }$ when $\deg \beta = \ell -1$. Next, choose $\beta$ of degree $\ell-2$ and $\alpha$ of degree 1. By the same
argument, the coefficient of $\alpha z x_j \beta$ vanishes, hence $\sum_{1\leq i \leq n} \lambda_{i,\alpha, \beta } f_{ji}=0$ for $j=1, \ldots , n$,
and we conclude again that all the elements $\lambda_{i,\alpha, \beta }$ vanish when $\deg \beta = \ell -2$. Continuing the
process, we arrive to $S=0$. One proves similarly that the right $A$-module $zA$ is free generated by $z$.
\Edm
\\

Lemma \ref{free1} has the following consequence: for any $a\in A$, there exists a
unique element $\sigma (a)\in A$ (resp. $\sigma'(a)\in A$) such that $za=\sigma (a)z$ (resp. $az= z\sigma'(a)$). Clearly, $\sigma$ and $\sigma'$ are
automorphisms of  the $k$-algebra $A$, inverse to each other. Now we are going to describe $B$ as a skew polynomial algebra with coefficients in $A$. It
suffices to do it for $T_A(Az) = \bigoplus_{p\geq 0} (Az)^{\otimes_A p}$. Proceeding by induction on $p$, we prove easily the following.
\Blm \label{free2}
The left (or right) $A$-module $(Az)^{\otimes_A p}$ is free generated by
$z^{\otimes_A p}=z\otimes_A \cdots \otimes_A z$ ($p$ times) for any $p\geq 1$.
\Elm

Consequently, any element of $T_A(Az)$ is uniquely written as a finitely supported
sum $\sum_{p\geq 0} a_p z^p$ where $a_i\in A$ and where we set $z^p= z^{\otimes_Ap}$. The product in $T_A(Az)$ is determined by the product in $A$
and by the relations $$z^p a= \sigma ^{p} (a) z^p,$$
for any $a\in A$ and $p\geq 1$. We have obtained the following proposition (on skew polynomial rings, see
e.g.~\cite{bg:aqg} p. 8-9).
\Bpo  \label{skew}
For any non-degenerate $f$, the $k$-algebra $B=B(f)$ is isomorphic to the skew
polynomial $k$-algebra $A[z;\sigma]$ defined over the $k$-algebra $A$ by $z$ and the $k$-automorphism $\sigma$ of $A$.
\Epo

The isomorphism $B \cong A[z;\sigma]$ is an isomorphism of graded algebras, knowing
that $A$ is graded and $z$ has degree 1. For any non-zero element $a=\sum_{p\geq 0} a_p z^p$ in $B$, its degree $\deg_z(a)$
in $z$ is the highest $p$ such that $a_p\neq 0$. One has $\deg_z(az)=\deg_z(za)=\deg_z(a)+1$, so that $z$ is not a zero-divisor in $B$.
Moreover it is easy to deduce the Hilbert series of $B$ from the Hilbert series of~$A$:
\begin{equation}
h_B(t)=\frac{h_A(t)}{1-t}.
\end{equation}
Using Proposition \ref{recall}, we get
\begin{equation} \label{hs1}
h_B(t)=(1-(n+1)t + (n+1)t^2 - t^3)^{-1} \ \ \mbox{if} \ n\geq 2,
\end{equation} \begin{equation} \label{hs2}
h_B(t)=(1-2t+2t^2-2t^3 + \cdots )^{-1} \ \ \mbox{if} \ n=1.
\end{equation}
\Bpo
Assume that $f$ is non-degenerate. The quadratic algebra $B=B(f)$ is Koszul. The global
dimension of $B$ is equal to 3 if $n\geq 2$, to $\infty$ if $n=1$. Moreover $B$ is a domain if and only if $n\geq 2$.
\Epo
\Bdm
The global dimension of the graded algebra $B$ is immediately derived from the
expression (\ref{hs1}) or (\ref{hs2}) of its Hilbert series. According to a standard result on Koszul algebras~\cite{bafr:kalg}, the Koszulity of $B$ comes from
the Koszulity of $A=B/(z)$ because the element $z$ of degree 1 in $B$ is normal and is not a zero-divisor (see also Example 1, p. 33 in~\cite{popo:quad}). The last statement is clear from $B \cong A[z;\sigma]$ and from 8) in Proposition \ref{recall}.
\Edm
\\

Assume that $f$ is non-degenerate as before, with $n\geq 2$. The Hilbert series (\ref{hs1})
of the Koszul algebra $B$ shows that $$\dim((R_B \otimes V_B)\cap (V_B\otimes R_B))=1.$$ Since $c(w)$ is an element of $(R_B \otimes V_B)\cap (V_B\otimes R_B)$ by the
non-commutative Euler relation (\ref{euler}), one has $$(R_B \otimes V_B)\cap (V_B\otimes R_B)= k\, c(w).$$
Therefore the bimodule Koszul complex of $B$ is the following complex $K_w$
\begin{equation} \label{resKo}
0 \longrightarrow B \otimes k\,c(w) \otimes B
\stackrel{d_3}\longrightarrow B \otimes R_B \otimes B \stackrel{d_2}\longrightarrow
B \otimes V_B \otimes B \stackrel{d_1}\longrightarrow B \otimes B
\end{equation} and $K_w \stackrel{\mu}\longrightarrow B$ is the Koszul resolution of $B$ where $\mu
: B \otimes B \rightarrow B$ is the multiplication. To simplify notation, set $x_{n+1}=z$ and $r_i=\partial_{x_i}(w)$ for $i=1, \ldots , n+1$.
Beside the cyclic partial derivative $\partial_{x_i}$ defined on $F_{cyc}$, there is the ``ordinary'' partial derivative $\frac{\partial}{\partial x_i} : F
\rightarrow F\otimes F$ (see~\cite{vdb:ober}) defined on any monomial $a$ by $$\frac{\partial a}{\partial x_i} = \sum_{a=ux_iv} u\otimes v,$$
which will be written as $$\frac{\partial a}{\partial x_i} = \sum_{1,2} \left(\frac{\partial a}{\partial x_i}\right)_1
\otimes \left(\frac{\partial a}{\partial x_i}\right)_2 \ .$$
Then the differential of the Koszul complex $K_w$ is given by
\begin{equation}
d_1(1\otimes x_i\otimes 1)=x_i\otimes 1 - 1\otimes x_i,
\end{equation}
\begin{equation}
d_2(1\otimes r_i \otimes 1)=\sum_{1\leq j \leq n+1} \sum_{1,2} \left(\frac{\partial r_i}{\partial x_j}\right)_1
\otimes x_j \otimes \left(\frac{\partial r_i}{\partial x_j}\right)_2 \end{equation}
\begin{equation}
d_3(1 \otimes c(w) \otimes 1)=\sum_{1\leq j \leq n+1} x_j\otimes r_j\otimes 1 - 1\otimes r_j \otimes x_j.
\end{equation}
The following theorem is an immediate consequence of Keller's theorem (Theorem 4.8 in~\cite{bk:compl}), since the algebra $B(f)$ is a 3-Calabi-Yau completion in sense of Keller (see 6.8 in~\cite{bk:compl}). The reader could be interested in a proof not using Keller's completions, so we present below such a proof based on elementary homological algebra. \Bte \label{calyau}
For any non-degenerate $f$ and for any $n\geq 2$, the algebra $B=B(f)$ is
3-Calabi-Yau. In other words, the potential $w=fz$ is 3-Calabi-Yau.
\Ete
\Bdm
It suffices to prove that the complex $K_w$ is self-dual with respect to the functor
$(\,\cdot\,)^{\vee}=Hom_{B-B}(\,\cdot\,, B\otimes B)$ (\cite{rbrt:pbw}, Lemma 3.7). Denote by $E^{\ast}$ the dual space of a $k$-linear
space $E$. It is easy to compute the complex $K_w^{\vee}$:
\begin{equation}
B \otimes B
\stackrel{d_1^{\ast}}\longrightarrow B \otimes V_B^{\ast} \otimes B
\stackrel{d_2^{\ast}}\longrightarrow
B \otimes R_B^{\ast} \otimes B \stackrel{d_3^{\ast}}\longrightarrow B \otimes
kc(w)^{\ast} \otimes B \longrightarrow 0
\end{equation}
where the differential is given on the dual basis by
\begin{equation}
d_1^{\ast}(1\otimes 1)=\sum_{1\leq j \leq n+1} x_j\otimes x_j^{\ast} \otimes 1 - 1\otimes
x_j^{\ast} \otimes x_j,
\end{equation}
\begin{equation}
d_2^{\ast}(1\otimes x_i^{\ast}\otimes 1)=\sum_{1\leq j \leq n+1} \sum_{1,2} \left(\frac{\partial
r_j}{\partial x_i}\right)_2 \otimes r_j^{\ast} \otimes \left(\frac{\partial r_j}{\partial x_i}\right)_1 \end{equation}
\begin{equation}
d_3^{\ast}(1 \otimes r_i^{\ast} \otimes 1)= x_i \otimes c(w)^{\ast}\otimes 1 - 1\otimes c(w)^{\ast}
\otimes x_i.
\end{equation}
Consider the diagram (in which we have omitted the symbols $\otimes$)
\begin{eqnarray} \label{self}
B(kc(w))B \stackrel{d_{3}}{\longrightarrow} & BR_{B} B
\stackrel{d_{2}}{\longrightarrow} & BV_{B} B \stackrel{d_{1}}{\longrightarrow}  Bk B
\nonumber \\
 f_{3} \downarrow \ \ \ \ \ \ \  & f_{2} \downarrow \ \ \ \ \ \ \ &  f_{1}
\downarrow \ \ \ \ \ \ \ \ \  f_{0} \downarrow \\
Bk B \stackrel{d_{1}^{\ast}}{\longrightarrow} & BV^{\ast}_{B} B
\stackrel{d_{2}^{\ast}}{\longrightarrow}
& BR^{\ast}_{B} B \stackrel{d_{3}^{\ast}}{\longrightarrow}  B(kc(w)^{\ast}) B
\nonumber \end{eqnarray}
where the $B$-bimodule isomorphisms $f_i$ are given by
\begin{equation}
f_0(1)=c(w)^{\ast},\ f_1(x_i)= r_i^{\ast}, \ f_2(r_i)=x_i^{\ast}, \ f_3(c(w))= 1.
\end{equation}
The diagram (\ref{self}) is \emph{commutative}. In fact, it is immediate to check
that the left square and the right square are commutative. Moreover the commutativity of the central square is a straightforward consequence of
the following non-commutative Hessian formula~\cite{vdb:ober}
\begin{equation}
\tau \left(\frac{\partial^2 w}{\partial x_i \partial x_j}\right)= \frac{\partial^2 w}{\partial
x_j \partial x_i}
\end{equation} where we set $\frac{\partial^2 }{\partial x_i \partial x_j}=
\frac{\partial}{\partial x_i} \circ \partial_{x_j}$ and where $\tau :F\otimes F \rightarrow F\otimes F$ is the flip. Thus $f:K_w \rightarrow
K_w^{\vee}$ is a complex isomorphism, so that $K_w$ is self-dual.
\Edm
\\

We complete the properties of $B=B(f)$ by the following.
\Bpo  \label{comple}
Assume that $f$ is non-degenerate and $n\geq 2$. \\ 1) $B$ is AS-Gorenstein.
\\ 2) If $f$ is alternating, then $z$ is central in $B$. The converse holds if the
characteristic of $k$ is $\neq 2$.
\\ 3) The Gelfand-Kirillov dimension GK.dim$(B)$
of $B$ is finite if and only if $n=2$, and in this case GK.dim$(B)=3$.
\\ 4) $B$ is left (or right) noetherian if and only if $n=2$.
\Epo
\Bdm \ 1) It is an immediate consequence of the Calabi-Yau property (\cite{rbrt:pbw},
Proposition 4.3). 2) The first claim comes from Example 2.4. For the converse, if $z$ is central in
$B$, then we have in $B$ the equalities $\sum_{1\leq j \leq n} (f_{ij}+ f_{ji}) x_j z = 0$ for $i=1,\ldots , n$. But the
$A$-module $Az$ is free, hence $x_1z, \ldots x_nz$ are linearly independent in $B$, and we conclude that $f_{ij}+f_{ji}=0$ for any $i$ and $j$.

3) For $n=2$, we have $h_B(t)=(1-t)^{-3}$, hence GK.dim$(B)=3$. If $n>2$, the
polynomial $1-(n+1)t + (n+1)t^2 -t^3$ has a real root between 0 and 1, hence GK.dim$(B)=\infty$ .

4) If $n=2$, the algebras $A$ are noetherian~\cite{jz:nnoeth}. Thus $B\cong A[z;\sigma]$ is noetherian by a standard result~\cite{bg:aqg}. If $n>2$, then GK.dim$(B)=\infty$, so we conclude that $B$ is not (left or right) noetherian by using a theorem due to
Stephenson and Zhang~\cite{sz:gro}. Let us note that the equivalence in 4) can be deduced from the analogous equivalence for $A$ (see 7) in Proposition \ref{recall}) and from a general result (\cite{lbsvdb:central}, Corollary 2.3).
\Edm

\setcounter{equation}{0}

\section{Classification of the algebras $B$}
The aim of this section is to classify the algebras $B=B(M)$ up to isomorphims of
graded algebras. Denote by $A(V,R)$ the quadratic graded algebra defined by a
space of generators $V$ and a space of relations $R$ (subspace of $V\otimes V$). We write
$A(V,R)\cong A(V',R')$ if there exists an isomorphism of graded algebras from
$A(V,R)$ to $A(V',R')$. It is equivalent to say that there exists a linear
isomorphism $\varphi: V \rightarrow V'$ such that $(\varphi \otimes \varphi )(R)=R'$ (this property is the definition of
isomorphisms between two quadratic graded algebras given in~\cite{man:kalg, bdvw:homog}). When
$M\in M_n(k)$, we begin to classify up to isomorphism the graded algebras $A(M)$ defined at the beginning of Section 2. Using the same notation, let $A(M)$ be the quadratic
algebra defined as the quotient of $k\langle x_1, \ldots , x_n\rangle$ by the two-sided ideal generated by $f$, where $f=\sum _{1\leq i,j \leq n} f_{ij}
x_i x_j$ and $M=(f_{ij})_{1\leq i,j \leq n}$. It is clear that for any matrices $M_1$ and $M_2$, the equality $A(M_1)=A(M_2)$ holds if and only if there exists a non-zero scalar $\nu$ such that $M_1=\nu M_2$. We have \begin{equation} \label{form1}
f= \left( \begin{array}{ccc}
x_1  & \ldots & x_n  \end{array} \right) M \left( \begin{array}{c}
x_1  \\ \vdots \\ x_n \end{array} \right).
\end{equation}
The isomorphism class of $A(M)$ is the class of the algebras $A(N)$ which can be defined by $M$ with respect to another basis $(x'_1, \ldots , x'_n)$. It means that $A(N)$ is the quotient $$k\langle x_1, \ldots , x_n\rangle/(f')$$
where \begin{equation} \label{form11}
f'= \left( \begin{array}{ccc}
x'_1  & \ldots & x'_n  \end{array} \right) M \left( \begin{array}{c}
x'_1  \\ \vdots \\ x'_n \end{array} \right),
\end{equation}
in which $f'$ has to be written in the variables $x_1, \ldots, x_n$.
Considering $P\in GL_n(k)$ such that $$\left( \begin{array}{c}
x'_1 \\ \vdots \\ x'_n \end{array} \right)
= P \left( \begin{array}{c}
x_1  \\ \vdots \\ x_n  \end{array} \right),$$
we obtain
\begin{equation} \label{form12}
f'= \left( \begin{array}{ccc}
x_1  & \ldots & x_n  \end{array} \right)\,^tPM P \left( \begin{array}{c}
x_1  \\ \vdots \\ x_n \end{array} \right).
\end{equation}
Thus the matrices $N$ are the matrices of the form $N=\nu \,^tPM P$ for any non-zero scalar $\nu$ and any $P\in GL_n(k)$. Two matrices $M$ and $N$ in
$M_n(k)$ are said to be \emph{congruent} if there exists $P\in GL_n(k)$ such that
$N=\,^tPM P$. Note that we can choose $\nu = 1$ in $N=\nu \,^tPM P$ if any element of $k$ is a square in $k$, but if $\nu$ is not a square in $k$, $\nu$ is not congruent to 1 (for $n=1$). We have obtained the following result (when $k=\mathbb{C}$, see \cite{rb:gera}, end of Section 5). \Bpo \label{classA}
Let $k$ be a field and $n\geq 1$. Let $M$ and $N$ be in $M_n(k)$. Then $A(M)\cong
A(N)$ if and only if $M$ is congruent to a non-zero scalar multiple of $N$. If any element of $k$ is a square in $k$, then $A(M)\cong A(N)$ if and only if $M$ and $N$ are congruent.
\Epo

Now we are interested in the classification of the graded algebras $B$ up to
isomorphism. We shall prove that this classification is the same as for the algebras
$A$. Keeping the algebra $A(M)$ as above, we define the quadratic algebra $B(M)$ by the potential $w=fz$ as at the beginning of Section 2. The isomorphism
class of $B(M)$ is the class of the algebras $B(N)$ which can be defined by $M$ with respect to another basis $(x'_1, \ldots , x'_n,z')$. To be more precise, let $\Lambda \in GL_{n+1}(k)$ be such that
\begin{equation} \label{Lambda}
\left( \begin{array}{c}
x'_1 \\ \vdots \\ x'_n  \\ z'
\end{array} \right) = \Lambda \left( \begin{array}{c}
x_1  \\ \vdots \\ x_n \\ z
\end{array} \right). \end{equation}
Note that $z$ can occur in $x'_1, \ldots , x'_n$ and that $x_1, \ldots , x_n$ can occur in $z'$. Write $f'$ given by (\ref{form11}) and $w'=f'z'$ in the variables $x_1, \ldots, x_n,\,z$. Consider the quadratic algebra $$\Lambda \cdot B(M):= k\langle x_1, \ldots , x_n, \, z\rangle/(\partial_{x'_1}(w'),
\ldots , \partial_{x'_n}(w'), f').$$
Then we have to find the algebras $B(N)$ which are equal to $\Lambda \cdot B(M)$ for some $\Lambda \in GL_{n+1}(k)$. Clearly the isomorphism class of $B(0)=k\langle x_1, \ldots , x_n, \, z\rangle$ is the singleton $\{B(0)\}$.
\Bpo \label{direct}
If $N=\nu \,^tPM P$ for some non-zero scalar $\nu$ and $P\in GL_n(k)$, then $B(N)=\Lambda \cdot B(M)$, where $\Lambda = \left( \begin{array}{cc}
P & 0 \\ 0 & \nu \end{array} \right).$
\Epo
\Bdm
Let $(x'_1, \ldots , x'_n,z')$ be as in (\ref{Lambda}) with $\Lambda = \left( \begin{array}{cc}
P & 0 \\ 0 & \nu \end{array} \right)$. According to (\ref{rel}), we have
\begin{equation} \label{relprime}
\left( \begin{array}{c}
\partial_{x'_1}(w') \\ \vdots \\ \partial_{x'_n}(w') \end{array} \right)
= \left( \begin{array}{cc}
M & \,^t M
\end{array} \right) \left( \begin{array}{c}
x'_1 z' \\ \vdots \\ z' x'_n \end{array} \right).
\end{equation}
Consequently, from $$\left( \begin{array}{c}
x'_1 \\ \vdots \\ x'_n  \end{array} \right) = P
\left( \begin{array}{c}
x_1  \\ \vdots \\ x_n \end{array} \right)$$
and $z' = \nu z$, we obtain \begin{equation} \label{rel2}
\,^tP \left( \begin{array}{c}
\partial_{x'_1}(w') \\ \vdots \\ \partial_{x'_n}(w') \end{array} \right)
= \left( \begin{array}{cc}
\, N & \,^tN
\end{array} \right) \left( \begin{array}{c}
x_1 z \\ \vdots \\ z x_n \end{array} \right),
\end{equation}
and \begin{equation} \label{form2}
f'= \left( \begin{array}{ccc}
x_1  & \ldots & x_n  \end{array} \right) \,^tPMP  \left( \begin{array}{c}x_1  \\ \vdots \\ x_n \end{array} \right).
\end{equation}
Thus the spaces
of relations of $\Lambda \cdot B(M)$ and $B(N)$ are equal.
\Edm \\

For any $\Lambda \in GL_{n+1}(k)$, write $\Lambda = \left( \begin{array}{cc} \widetilde{\Lambda} & c \\ \ell & \lambda \end{array} \right)$ where $\widetilde{\Lambda} \in M_n(k)$. In order to illustrate the proof of Theorem \ref{classB2} below when $\widetilde{\Lambda}$ is not invertible, the next proposition provides some explicit examples.
\Bpo \label{ntexample}
Let $\Lambda$ be as in (\ref{Lambda}) with the basis change $x'_1= z, \ x'_2=x_2,\ \ldots ,\ x'_n=x_n, \newline z'=x_1$. For any matrix \begin{equation} \label{specialM}
M = \left( \begin{array}{cccc} 0 & f_{12} & \ldots & f_{1n} \\
f_{21} & 0 & \ldots & 0 \\
\vdots & \vdots & \vdots & \vdots \\
f_{n1} & 0 & \ldots & 0
\end{array} \right),
\end{equation}
we have $\Lambda \cdot B(M)=B(^tM)$.
\Epo
\Bdm
Performing the basis change in $f'$ given by (\ref{form11}) and in $w'=f'z'$, we obtain
$f'= \sum_{2\leq j \leq n}(f_{1j}zx_j + f_{j1}x_jz)$ and $w'= \sum_{2\leq j \leq n}(f_{1j}zx_jx_1 + f_{j1}x_jzx_1)$. Therefore the relations of $\Lambda \cdot B(M)$ are
$$\partial_{x'_1}(w')=\partial_{z}(w')=\sum_{2\leq j \leq
n}(f_{1j}x_jx_1 + f_{j1}x_1x_j),$$
$$\partial_{x'_i}(w')=\partial_{x_i}(w')= f_{1i}x_1z + f_{i1}zx_1 ,\ \ 2 \leq i \leq n,$$
$$f'= \sum_{2\leq j \leq n}(f_{1j}zx_j + f_{j1}x_jz).$$
Since $B(^tM)$ is defined by the potential $gz= \sum_{2\leq j\leq n}(f_{1j}x_jx_1z + f_{j1}x_1x_jz)$, the relations of $B(^tM)$ are  $$\partial_{x_1}(gz)= \sum_{2\leq j \leq n}(f_{1j}zx_j + f_{j1}x_jz),$$$$\partial_{x_i}(gz)= f_{1i}x_1z + f_{i1}zx_1 ,\ \ 2 \leq i \leq n,$$
$$g=\sum_{2\leq j \leq n}(f_{1j}x_jx_1 + f_{j1}x_1x_j).$$
Thus the spaces of relations of $\Lambda \cdot B(M)$ and $B(^tM)$ are equal. \Edm
\\

Actually, for the matrix $\Lambda$ of Proposition \ref{ntexample}, one can show that any matrix $M$ such that there exists $N$ satisfying $\Lambda \cdot B(M)=B(N)$ is of the form (\ref{specialM}), and then $N$ is a non-zero scalar multiple of $^tM$. For the proof of the classification theorem of the algebras $B(M)$ (Theorem \ref{classB2} below), one could probably do things intrinsically by using matrices $\Lambda$, but following the suggestion of the referee, we prefer a straightforward approach which is structurally simpler.
\Bte \label{classB2}
Let $k$ be a field and $n\geq 1$. Let $M=(f_{i,j})$ and $N=(f'_{i,j})$ be in $M_n(k)$. Then $B(M)\cong
B(N)$ if and only if $M$ is congruent to a non-zero scalar multiple of $N$. If any element of $k$ is a square in $k$, then $B(M)\cong B(N)$ if and only if $M$ and $N$ are congruent. \Ete
\Bdm
According to Proposition \ref{direct}, it remains to prove the necessity of the condition when no assumption is made on $k$. Assume that there exists an isomorphism of graded algebras $\phi : B(M) \rightarrow B(N)$. As noted just before Proposition \ref{direct}, $M=0$ if and only if $N=0$, so that we assume that $M$ and $N$ are non-zero. Summing implicitly over repeated indices, we have $\phi(x_i)=\alpha_{i,j} x_j + \beta_i z$ and $\phi(z)=\gamma _jx_j + \delta z$ for some well-determined scalars $\alpha_{i,j}$, $\beta_i$, $\gamma _j$ and $\delta$. Denote by $\bar{\phi}: B(N) \rightarrow B(M)$ the map inverse of $\phi$, and let $\bar{\alpha}_{i,j}$, $\bar{\beta_i}$, $\bar{\gamma}_j$, $\bar{\delta}$ be the corresponding scalars for the map $\bar{\phi}$. As $\bar{\phi}\phi = 1_{B(M)}$, we have \begin{equation} \label{rcoef1}
\begin{array}{ll}
\alpha_{i,j}\bar{\alpha}_{j,k} + \beta_i\bar{\gamma}_k= \delta _{i,k},  & \alpha_{i,j}\bar{\beta}_j + \beta_i \bar{\delta}= 0,
\\
\gamma_i \bar{\alpha}_{i,j} + \delta \bar{\gamma}_j= 0,  \ \ \ \ \ & \gamma_i \bar{\beta}_i + \delta \bar{\delta}= 1; \end{array}
\end{equation}
of course, there are similar equations coming from $\phi \bar{\phi}=1_{B(N)}$.

As $f_{i,j}x_ix_j=0$ in $B(M)$, we have $f_{i,j}(\alpha_{i,k}x_k + \beta_iz)(\alpha_{j,l}x_l + \beta _jz)=0$ in $B(N)$. Since the algebra $B(N)$ is graded by $\deg_z$, looking at the homogeneous components of degrees zero and two of this equality tells us that in $B(N)$ \begin{equation*} \begin{array}{ll}
f_{i,j} \alpha_{i,k} \alpha_{j,l} x_k x_l = 0,  & f_{i,j} \beta_i \beta_j z^2 = 0. \end{array}
\end{equation*}
The first of these tells us that there is a scalar $\lambda$ such that \begin{equation} \label{rcoef2}
f_{i,j} \alpha_{i,k} \alpha_{j,l} = \lambda f'_{k,l}
\end{equation}
and the second one $-$ since $z^2\neq0$ in $B(N)$ $-$ that
\begin{equation} \label{rcoef3}
f_{i,j} \beta_i \beta_j = 0. \end{equation}

If the matrix $(\alpha_{i,j})$ is invertible, then $\lambda$ cannot be zero, for we would then have that $M=0$ in view of equation (\ref{rcoef2}), and that same equation tells us that $M$ and $\lambda N$ are congruent. We may therefore assume that it is not and, by symmetry, that neither is the matrix $(\bar{\alpha}_{i,j})$. This assumption implies that there exist scalars $\zeta_i$, not all zero, such that $\zeta_i \alpha_{i,j}=0$ for all $j$, and then
$$\zeta_j x_j= \bar{\phi} \phi (\zeta_i x_i)= \bar{\phi}(\zeta_i\alpha_{i,j} x_j + \zeta_i \beta_i z) = \zeta_i \beta_i (\bar{\gamma}_j x_j + \bar{\delta}z).$$
It follows that $\zeta_i \beta_i \bar{\gamma}_j = \zeta _j$ for all $j$ and that $\zeta_i \beta_i \bar{\delta}=0$; since not of all of the $\zeta_j$ are zero, this implies that $\zeta_i \beta_i \neq 0$ and therefore $\bar{\delta}=0$. By symmetry, we also have $\delta =0$. As $f_{i,j} x_jz + f_{j,i} zx_j$ is zero in $B(M)$ and the homogeneous component of degree zero of its image under $\phi$ is $(f_{i,j} \alpha_{j,k} \gamma_l + f_{j,i} \alpha_{j,l} \gamma_k)x_k x_l$, this must be a multiple of $f'_{k,l} x_kx_l$, and we see that there exist scalars $\mu_i$ such that \begin{equation} \label{rcoef4}
f_{i,j} \alpha_{j,k} \gamma_l + f_{j,i} \alpha_{j,l} \gamma_k= \mu_i f'_{k,l}. \end{equation}

If we multiply (\ref{rcoef2}) by $\bar{\beta}_k$ and sum over $k$, we see that $\lambda f'_{k,l} \bar{\beta}_k= f_{i,j} \alpha_{i,k} \alpha_{j,l} \bar{\beta}_k=0$, because of the second equation in (\ref{rcoef1}). If we had $\lambda \neq 0$, so that $f'_{k,l} \bar{\beta}_k=0$, we could multiply equation (\ref{rcoef4}) by $\bar{\beta}_k$, sum over $k$, and $-$ using that $f'_{k,l} \bar{\beta}_k=0$ and $\gamma_k \bar{\beta}_k=1$ $-$ we would have that $f_{j,i} \alpha_{j,l}=0$; equation (\ref{rcoef2}) would then tell us that in fact $f'=0$, against our hypothesis. Therefore $\lambda=0$. Let $\omega_{i,j}=\alpha_{i,j} + \beta_i\gamma_j$ and similarly for $\bar{\omega}_{i,j}$. Using equations (\ref{rcoef1}) and the fact that $\delta =\bar{\delta}=0$, we immediately see that the matrix $(\omega_{i,j})$ is invertible with inverse $(\bar{\omega}_{i,j})$. Also we have
\begin{equation} \label{rcoef5}
f_{i,j} \omega_{i,k} \omega_{j,l} = f_{i,j} \alpha_{i,k} \alpha_{j,l} + f_{i,j} (\alpha_{i,k}\beta _j\gamma_l + \alpha_{j,l}\beta_i\gamma_k) +  f_{i,j}\beta_i\beta_j\gamma_k\gamma_l. \end{equation}
The first term is zero because of (\ref{rcoef2}) and the third one vanishes according to (\ref{rcoef3}), so that we have
$$f_{i,j} \omega_{i,k} \omega_{j,l} = f_{i,j}\alpha_{i,k}\beta _j\gamma_l + f_{i,j}\alpha_{j,l}\beta_i\gamma_k = (f_{j,i}\alpha_{j,k}\gamma_l + f_{i,j}\alpha_{j,l}\gamma_k) \beta_i = \mu_i \beta_i f'_{l,k},$$
in view of equation (\ref{rcoef4}). Since the matrix $(\omega_{i,j})$ is invertible and not all the $f_{i,j}$ are zero, this implies that $\nu:=\mu_i \beta_i\neq 0$ and we conclude that $M$ and the transpose of $\nu N$ are congruent. Since any matrix is congruent to its transpose~\cite{doik:congr}, it follows that $M$ is congruent to $\nu N$ itself.
\Edm  
\setcounter{equation}{0}

\section{Special properties of $B(M)$ when $n=2$}
Throughout the sequel of this paper, we assume that $k=\mathbb{C}$ (actually, it
would be sufficient to assume that $k$ is algebraically closed of characteristic zero). An explicit parameter space for the $GL(2,\mathbb{C})$-orbits of congruent matrices in $GL(2, \mathbb{C})$ is given by Dubois-Violette
in~\cite{dv:multi} (for the classification of matrices under congruence in general, the reader may consult~\cite{cw:bili, ser:classif}). According to the previous section, this parameter space forms a ``moduli space'' for the algebras $B(M)$ when $M$ runs over $GL(2, \mathbb{C})$.
Recall the description in three types of this parameter space as stated in (\cite{dv:multi}, end of Section 2). In each type, we give the potential $w$ and
the relations of $B(M)$. The generators are denoted by $x$, $y$ and $z$. The types depend on the rank \textbf{rk} of the
symmetric part $s(M)=\frac{1}{2} (M+\,^tM)$ of $M$. It turns out that the corresponding algebras $A(M)$ are exactly the AS-regular algebras of global
dimension 2~\cite{as:regular}. Therefore the three types classify the AS-regular algebras of global dimension 2 as well, and we have kept the samestandard terminology for the classification of our algebras $B(M)$.
\\ \\
\emph{First type (classical type)}: $\textbf{rk}=0$. There is only one orbit, which
is the orbit of $M=\left( \begin{array}{cc} 0 & -1 \\ 1 & 0 \end{array} \right)$. One has $w=(yx-xy)z$. The relations of $B(M)$ are the
following: $zy=yz$, $xz=zx$, $yx=xy$. Hence $B(M)=\mathbb{C}[x,y,z]$ the commutative polynomial algebra. Note that it is the unique orbit such that $z$
is central in $B(M)$ (Proposition \ref{comple}), or such that $A(M)$ is Calabi-Yau (Proposition \ref{recall}).
\\ \\
\emph{Second type (Jordan type)}: $\textbf{rk}=1$. There is only one orbit, which is
the orbit of $M=\left( \begin{array}{cc} -1 & -1 \\ 1 & 0 \end{array} \right)$. One has $w=(yx-xy-x^2)z$. The relations of $B(M)$ are the
following: $zy=yz+2xz$, $xz=zx$, $yx=xy+x^2$.
\\ \\
\emph{Third type (quantum type)}: $\textbf{rk}=2$. The orbits are parametrized by
the set $\{q \in \mathbb{C} \setminus\{0,1\}\}/(q \sim q^{-1})$. These are the orbits of $M=\left( \begin{array}{cc} 0 & -q^{-1} \\ 1 & 0 \end{array}
\right)$. One has $w=(yx-q^{-1}xy)z$. The relations of $B(M)$ are the following: $xy=q yx$, $yz=q zy$, $zx= qxz$. So $B(M)$ is a quantum affine space.
\\
\Bpo\label{prop:basisofB}
For any $M\in GL(2, \mathbb{C})$, the algebra $B(M)$ is a quadratic AS-regular
algebra of global dimension 3, of type A non-generic (following the classification of Artin and Schelter~\cite{as:regular}). For any $M$ in the list of the above
classification, $(x^iy^jz^k)_{i\geq 0, j \geq 0, k \geq 0}$ is a basis of the vector space $B(M)$. \Epo
\Bdm
According to Proposition \ref{comple}, $B(M)$ is AS-Gorenstein with a finite
Gelfand-Kirillov dimension, hence is AS-regular. Since $B(M)$ is Calabi-Yau (Theorem \ref{calyau}), $B(M)$ is of type A (\cite{rbrt:pbw},
Proposition 5.4). Moreover, $B(M)$ is a skew polynomial ring over $A(M)$ in $z$ (Proposition \ref{skew}) and it is clear from each relation $f$ listed above
that $A(M)$ is a skew polynomial ring in $x$, $y$. Therefore $B(M)$ is a skew polynomial ring in $x$, $y$, $z$ (hence the basis of monomials
$x^iy^jz^k$), and the invariant $j$ of $B(M)$ is infinite (\cite{as:regular},
Theorem 6.11).  \Edm
\\ Among the AS-regular algebras of type A, the generic ones correspond to a finite
invariant $j$ and are called Sklyanin algebras~\cite{vdb:Kth}. Recall that Sklyanin algebras are defined
by three generators $x$, $y$, $z$, and three relations $\alpha yz+\beta zy+\gamma
x^2=0$, $\alpha zx+\beta xz+\gamma y^2=0$, $\alpha xy+\beta yx+\gamma z^2=0$, where $\gamma \neq 0$.
So the algebras $B(M)$ of the first or third type can be considered as limits of Sklyanin algebras by vanishing the parameter $\gamma$, but such a process is not possible for the second type.

In the classification of their quadratic regular algebras of global dimension 3,
Artin and Schelter define the invariant $j$ as the invariant of a certain cubic curve $\mathcal{C}$ in $\mathbb{P}^2$. As we shall see, the equation
$\phi =0$ of $\mathcal{C}$ is easily deduced from the potential $w$ of $B(M)$. More importantly, we shall interpret $\phi$ as a Poisson potential whose
associated Poisson bracket (defined as usual by $\nabla \phi$) is exactly the semi-classical limit of $B(M)$ viewed as a deformation of $\mathbb{C}[x,y,z]$.
Let us begin by the description of the curve $\mathcal{C}$ as in~\cite{as:regular}.  To avoid confusion with our notation, let us replace the notation $M$, $w$, $Q$ used
in~\cite{as:regular} by $\mathcal{M}$, $\mathcal{W}$, $\mathcal{Q}$.
The $3\times 3$ matrix $\mathcal{M}$ is defined by \begin{equation} \label{matm}
\left( \begin{array}{c}
\partial_{x}(w) \\ \partial_{y}(w) \\ f
\end{array} \right)
= \mathcal{M} \left( \begin{array}{c}
x \\ y \\ z \end{array} \right),
\end{equation} so that we have $$\mathcal{M}=\left( \begin{array}{cc} \,^tM z & M \left(
\begin{array}{c} x \\ y \end{array} \right)
\\ \left( \begin{array}{cc} x & y \end{array} \right)M & 0 \end{array} \right).$$ The element $\mathcal{W}$ is defined by $\mathcal{W}= X \mathcal{M} \,^tX$, where
$X=\left( \begin{array}{ccc} x & y & z \end{array} \right)$. From (\ref{matm}) and the non-commutative Euler relation (\ref{euler}), we deduce that \begin{equation} \label{aspot}
\mathcal{W}= c(w),
\end{equation}
that is, $\mathcal{W}$ is precisely the potential $\overline{w}$ once identified to
its cyclic sum. The matrix $\mathcal{Q}$ is defined by $X \mathcal{M} = \left( \begin{array}{ccc} \partial_{x}(w) & \partial_{y}(w) & f \end{array}
\right)\,^t\mathcal{Q}$. Then (\ref{euler}) implies that $\mathcal{Q}$ is equal to the identity matrix and we recover the fact that $B(M)$ is of type A (by
definition). Next, define the subvariety $\mathcal{C}$ of $\mathbb{P}^2$ by its equation
$S(\mathcal{W})=0$, where $S(\mathcal{W})$ denotes the symmetrization of the element $\mathcal{W}$. The symmetrization consists
in replacing products of variables in the free algebra by the same  products in the polynomial algebra. In our situation, $S(\mathcal{W})=3S(w)=3S(f)z$.
Set $\phi =S(f)z$. The above classification in three types allows us to limit
ourselves to the following matrices \begin{equation} \label{matM}
M=\left( \begin{array}{cc} a & b \\ 1 & 0 \end{array} \right), \ b\neq 0, \end{equation}
so that \begin{equation} \label{potpois}
w= (ax^2+bxy + yx)z, \ \ \phi= (ax^2 + (b+1)xy)z.
\end{equation}
In the classical type, one has $\mathcal{C}=\mathbb{P}^2$, and otherwise
$\mathcal{C}$ is the union of three straight lines, two of which coincide in the
Jordan type. \Bdf
The polynomial $\phi=\phi(M)$ is called the Poisson potential associated to the
algebra $B(M)$.
\Edf
This definition will be more natural when $B(M)$ will be viewed as a Gerstenhaber
deformation whose Poisson bracket $\{ \cdot ,\cdot \}$ will be defined on $\mathbb{C}[x,y,z]$ by the formulas $$\{x,y \}=\frac{\partial \phi}{\partial z}, \quad \{y,z \}=\frac{\partial
\phi}{\partial x},\quad \{z,x \}=\frac{\partial \phi}{\partial y},$$
i.e., $\{\cdot, \cdot\} = \frac{\partial \phi}{\partial z}\frac{\partial }{\partial
x}\wedge \frac{\partial }{\partial y}+\frac{\partial \phi}{\partial x}\frac{\partial
}{\partial y}\wedge \frac{\partial }{\partial z}+\frac{\partial \phi}{\partial
y}\frac{\partial }{\partial z}\wedge \frac{\partial }{\partial x}$.
\Bdf
The so-defined Poisson bracket $\{ \cdot ,\cdot \}$ is called the Poisson bracket
derived from the Poisson potential $\phi$.
\Edf
It is clear that $\phi$ belongs to the Poisson center of $\{ \cdot ,\cdot \}$
(because $\{x,\phi\}=\{y, \phi\}=\{z, \phi\}=0$), i.e., $\phi$ is a Casimir element
of the Poisson bracket. In the classical type, $\{ \cdot ,\cdot \}=0$. For the other
types, the next proposition shows that the Casimir element $\phi$ lifts to a
non-zero element $\Phi$ of the center of $B(M)$. Note that the potential $w$ lifts
$\phi$ in the free algebra, but $w$ vanishes in $B$ (it is a difference with
Sklyanin algebras for which the potential $w$ does not vanish in the algebra). Note
that $c(w)$ always vanishes in the algebra according to the non-commutative Euler
relation (\ref{euler}).

Most calculations in the algebra $B(M)$ where $M$ is given by (\ref{matM}) are based on the following relations \begin{equation} \label{relB}
\begin{array}{lll}
zy=(-a+\frac{a}{b})xz - byz,  & zx= -\frac{1}{b} xz, & yx= -a x^2 - b xy.
 \end{array}
\end{equation}
It is easy to verify that these relations are confluent. Therefore, Bergman's diamond lemma provides the basis $(x^iy^jz^k)_{i\geq 0, j \geq 0, k \geq 0}$ of $B$, as already stated in Proposition \ref{prop:basisofB}.   \Bpo \label{center}
The element $\Phi$ of $B(M)$ defined by $$\Phi=(ax^2 + (b+1) xy)z$$
belongs to the center of $B(M)$.
\Epo
\Bdm
Relations (\ref{relB}) allow us to decompose $x \Phi$ and $\Phi x$ in thebasis $(x^iy^jz^k)_{i\geq 0, j \geq 0, k \geq 0}$. The computations are
straightforward and they show that $x \Phi= ax^3z + (b+1) x^2yz=\Phi x$. In the same manner, we get $y \Phi =
a^2(b-1)x^3 z + a(b^2 - b-1)x^2yz - b(b+1) xy^2z=\Phi y$ and $z \Phi=ax^2z^2+(b+1)xyz^2=\Phi z$.
\Edm
\\

We are now interested in the comparison between the Hochschild homology of $B(M)$
and the Poisson homology of $(\mathbb{C}[x,y,z],\{ \cdot ,\cdot \})$. In the
classical type, the Hochschild homology of $\mathbb{C}[x,y,z]$ is well-known and
coincides with the Poisson homology for the corresponding $\{ \cdot ,\cdot \}$
(which is zero in this case). As proved in the next proposition, it is the
same for the quantum type if $q$ is not a root of unity. We are grateful to Wambst
for the explicit formulas of Hochschild homology. \Bpo \label{quantum}
Let $q$ be a non-zero complex number which is not a root of unity. Let $B$ be the
$\mathbb{C}$-algebra defined by the potential $w=(yx-q^{-1}xy)z$, i.e. defined by
the generators $x$, $y$, $z$, and the following relations $$xy=qyx, \ yz=qzy, \ zx=qxz.$$ Let $S=\mathbb{C}[x,y,z]$ be the polynomial $\mathbb{C}$-algebra in $x$, $y$, $z$,
endowed with the Poisson bracket derived from the Poisson potential $\phi=(1-q^{-1})xyz$. Denote by $Z=\mathbb{C}[xyz]$ the subalgebra of $S$ generated by
$xyz$.
Then the Hochschild homology of $B$ is isomorphic to the Poisson homology of $S$ and
is given by \\ \\
$HH_0 (B) \cong Z \oplus x\mathbb{C}[x] \oplus y \mathbb{C}[y] \oplus
z\mathbb{C}[z], $\\
$HH_1 (B) \cong ((\mathbb{C}[x]\oplus yzZ)\otimes x)\oplus ((\mathbb{C}[y]\oplus
xzZ) \otimes y) \oplus ((\mathbb{C}[z]\oplus xyZ)\otimes z),$\\
$HH_2(B) \cong (xZ\otimes (y\wedge z)) \oplus (yZ\otimes (z\wedge x))\oplus
(zZ\otimes (x\wedge y)),$\\
$HH_3(B) \cong Z \otimes (x\wedge y \wedge z),$\\
$HH_p(B) =0$ for any $p\geq 4$.
\Epo
\Bdm
On one hand, setting $x_1=x$, $x_2=y$, $x_3=z$, $q_{12}=q_{23}=q_{31}=q$ and
 applying Th\'eor\`eme 6.1 in \cite{wam:kosz}, we obtain the result for $HH_p(B)$.
diagonalizable in sense of Monnier, and the Poisson cohomology is computed in
\cite{mon:quad}.  Moreover, the duality $HP^{\bullet}(S) \cong HP_{3-\bullet}(S)$
holds \emph{a priori}, because a Poisson bracket deriving from a potential is
always unimodular, i.e., its modular class vanishes (see the definition of the modular vector field (and its class) in
\cite{wein:mod}, which is called the curl vector field in \cite{dufhar:curl}). The example 2.4 (i) of \cite{mon:quad}, then gives bases of the Poisson homology vector spaces of $S$ and the fact that they are isomorphic to the
Hochschild homology vector spaces of $B$. \Edm
\\
%

%

We have seen that the duality $HP^{\bullet}(S) \cong HP_{3-\bullet}(S)$ holds since the Poisson bracket derives from a
potential.
Since $B$ is 3-Calabi-Yau (Theorem \ref{calyau}), one also has the duality
$HH^{\bullet}(B) \cong HH_{3-\bullet}(B)$.    \setcounter{equation}{0}

\section{Poisson homology and Hochschild homology for the second type}

In this section, the field $k$ is $\mathbb{C}$ and $n=2$. The aim of this section is
to prove Theorem \ref{calc} of the Introduction, i.e., to determine the Hochschild
homology of the algebra $B$ of the second type (Jordan type):
$$
B = \mathbb{C}\langle x, y, z\rangle /(zy = yz+2xz, zx=xz, yx = x^2+xy).
$$
Let us recall the notation: $f=x^2+xy-yx$, $w=fz$, $r_1= \partial_{x}(w)$, $r_2=
\partial_{y}(w)$ and $r_3= \partial_{z}(w)$. We have $$r_1=xz+zx+yz-zy, \ r_2 = zx-xz, \ r_3 = x^2 + xy -yx.$$
In order to obtain the Hochschild homology of $B$, we will first see $B$ as the
deformation of a Poisson algebra for which we will determine the Poisson homology.
To do this, let us consider the filtration $F$ of $B$ given by the degree of $y$. In
other words, for this filtration, the degree of $y$ is 1 while the degrees of $x$
and $z$ are 0. It is clear from the monomial basis of $B$ (Proposition \ref{prop:basisofB}) that $gr_F(B)\simeq \mathbb
C[x,y,z]$, so the filtered algebra $B$ is almost commutative~\cite{kas:homcy}.
Moreover, $gr_F(B)$ is equipped with the Poisson bracket  defined by:
\begin{equation}
\lbrace z,y\rbrace =  2xz,\quad \lbrace z,x\rbrace =  0, \quad \lbrace y,x\rbrace = x^2,
\end{equation}
which is the Poisson structure derived from the Poisson potential $\phi=-x^2z$. In
the sequel, we will denote this Poisson algebra by $T=(\mathbb{C} \lbrack
x,y,z\rbrack, \lbrace\, ,æ\, \rbrace)$.

From the Koszul resolution $K_w$ (\ref{resKo}), we easily get the following Koszul
complex $B\otimes_{B^e}K_w$ associated to $B$ (where we have omitted the symbols
$\otimes$):
\begin{equation*}
0 \longrightarrow B(\mathbb{C}c(w))\stackrel{\tilde d_3}{\longrightarrow}  BR_B\stackrel{\tilde d_2}{\longrightarrow} BV_B \stackrel{\tilde d_1}{\longrightarrow} B \longrightarrow 0
\end{equation*}
where $c(w)=x^2z+xzx+zx^2 +xyz + yzx+zxy-yxz-xzy-zyx$ and where the differentials
are given, for $a\in B$, by:
\begin{equation*}
\tilde d_1 (a\otimes x) = \lbrack a, x\rbrack = ax - xa,\quad  \tilde d_1 (a\otimes y) = \lbrack a, y\rbrack,\quad  \tilde d_1 (a\otimes z) =
\lbrack a, z\rbrack, \end{equation*}
while
\begin{eqnarray*}
\tilde d_2 (a\otimes r_1) &=& (za+az)\otimes x + (za-az)\otimes y +
(ax+xa+ay-ya)\otimes z,\\
\tilde d_2 (a\otimes r_2) &=& (az-za)\otimes x + (xa-ax)\otimes z,\\
\tilde d_2 (a\otimes r_3) &=& (xa+ax+ya-ay)\otimes x + (ax-xa)\otimes y,
\end{eqnarray*}
and %
\begin{equation*}
\tilde d_3(a\otimes c(w)) = \lbrack a, x\rbrack \otimes r_1 + \lbrack a, y\rbrack
\otimes r_2 + \lbrack a, z\rbrack \otimes r_3.
\end{equation*}
This complex computes the Hochschild homology $HH_\bullet(B)$ of $B$. Notice that one can naturally identify the spaces of chains of this complex with $B$
or $B^3$. Indeed, one can write %
\begin{equation}\label{identHH}
B(\mathbb{C}c(w))\simeq B,\quad \hbox{ while }\quad BR_B\simeq B^3\; \hbox{ and }\;
BV_B\simeq B^3,
\end{equation}
 by identifying, for any $a_1, a_2, a_3\in B$, the element $a_1\otimes
r_1+a_2\otimes r_2+a_3\otimes r_3\in BR_B$ or the element  $a_1\otimes x+a_2\otimes
y+a_3\otimes z\in BV_B$ with $(a_1, a_2, a_3)\in B^3$.

As in $B$, we assume that in $V_B$, $R_B$ and $c(w)$, $y$ has degree 1 while $x$ and
$z$ have degree 0, so that we can consider the total filtration degree in $B\otimes_{B^e}K_w$. In Proposition \ref{po:gr} below, we show that $\tilde d_1$,
$\tilde d_2$ and $\tilde d_3$ have degree $-1$ for the total filtration degree. So
$B\otimes_{B^e}K_w$ will be a filtered complex (with differentials preserving the
degree) for the following shifts on the total filtration degree:
\begin{equation}\label{shift}
0 \longrightarrow B(\mathbb{C}c(w))(3)\stackrel{\tilde d_3}{\longrightarrow}  BR_B(2)\stackrel{\tilde d_2}{\longrightarrow} BV_B(1) \stackrel{\tilde d_1}{\longrightarrow} B \longrightarrow 0
\end{equation}
whose filtration is still denoted by $F$. %
Recall now the definition of the Poisson homology complex of  the Poisson algebra $T$.
First, the space of Poisson $1$-chains is the $T$-module of K\"ahler differentials
of $T$, denoted by $\Omega^1(T)$ and generated, as an $T$-module by the three
elements $d x$, $d y$, $d z$. Then, for $p\in\mathbb N^*$, the $T$-module of
K\"ahler $p$-differentials is the space of Poisson $p$-chains and is given by $\Omega^p(T)=\bigwedge^p\Omega^1(T)$. Of course, one has $\Omega^p(T)=\lbrace
0\rbrace$, as soon as $p\geq 4$. Notice that an
element $F_1 dx+F_2 dy+ F_3 dz\in \Omega^1(T)$ can naturally be identified with the element $(F_1, F_2, F_3)\in
T^3$. Similarly, an element $F_1 dy\wedge dz+F_2 dz\wedge dx + F_3 dx\wedge dy\in
\Omega^2(T)$ is identified with the element  $(F_1, F_2, F_3)\in T^3$, and an
element $F dx\wedge dy\wedge dz\in \Omega^3(T)$ with the element  $F\in T$. (Notice that in the following we will write elements in $T^3$ both as row vectors or as column vectors, depending on the context.)

Next, the Poisson boundary operator, $\delta_p:\Omega^p(T)\rightarrow\Omega^{p-1}(T)$,
called the Brylinski or Koszul differential, is given, for $F_0, F_1, \dots, F_p\in
T$, by (see \cite{bry:poi}):
\begin{eqnarray*}
\lefteqn{\delta_p(F_0\, d F_1\wedge\cdots\wedge d F_p)=
\sum_{i=1}^p(-1)^{i+1}\lbrace F_0,F_i\rbrace\,     d F_1\wedge\cdots\wedge\widehat{d F_i}\wedge\dots\wedge d F_p}\\
  &+&\displaystyle\sum_{1\leq i<j\leq p} (-1)^{i+j} F_0\, d\lbrace F_i, F_j\rbrace
        \wedge d F_1\wedge\cdots\wedge\widehat{d F_i}\wedge         \cdots\wedge\widehat{d F_j}\wedge\cdots\wedge d F_p,
\end{eqnarray*}%
where the symbol $\widehat{d F_i}$ means that we omit the term
$d F_i$.

Using the previous identifications, the Poisson homology complex of the Poisson
algebra $T$ can be written as:
\begin{equation*}
\setlength{\dgARROWLENGTH}{0.8cm}
\begin{diagram}
\node{0}\arrow{e,t}{}\node{T}\arrow{e,t}{\delta_3}\node{T^3}\arrow{e,t}{\delta_2}
              \node{T^3}\arrow{e,t}{\delta_1}\node{T}\arrow{e,t}{} \node{0}
\end{diagram}
\end{equation*}
with the differentials given, for $F\in T$, and $\vec F := (F_1, F_2, F_3)\in T^3$, by:
\begin{eqnarray}\label{eq:delta}
\delta_1 (\vec F) &=& \nabla\phi\cdot(\nabla\times \vec F) = \mathrm{Div}(\vec
F\times \nabla \phi),\nonumber\\
\delta_2(\vec F) &=& -\nabla(\vec F\cdot \nabla\phi)+\mathrm{Div}(\vec F)\nabla \phi,\\
\delta_3(F)&=& - \nabla F\times\nabla\phi\nonumber,
\end{eqnarray}
where $\mathrm{Div}(\vec F)=\frac{\partial F_1}{\partial x}+ \frac{\partial
F_2}{\partial y}+\frac{\partial F_3}{\partial z}\in T$, while, for all element $G\in
T$, we write $\nabla G := (\frac{\partial G}{\partial x}, \frac{\partial G}{\partial
y}, \frac{\partial G}{\partial z})\in T^3$. We also have respectively denoted by
$\cdot$, $\times$ and $\nabla\times$, the usual inner and cross products in $T^3$
and the curl operator.
Note that $\nabla \phi= (-2xz, 0, -x^2)$.

In proposition \ref{po:gr}, we show that
$gr_F(B\otimes_{B^e}K_w)$ is identified to the Poisson homology complex of the
Poisson algebra $T$. The next lemma will be useful.
\Blm\label{lm:form}
For all $k\in \mathbb N$, the following identities hold in $B$:
\begin{equation}\label{eq:yxk-zky}
yx^k = x^ky+ kx^{k+1}, \quad z^k y = yz^k+2kxz^k, \end{equation}
while
\begin{equation}\label{eq:ykx}
 y^kx= \sum_{j=0}^k \frac{k!}{j!}x^{k-j+1}y^j,
\end{equation}
and
\begin{equation}\label{eq:zyk}
zy^k = \sum_{\ell =0}^k (k-\ell+1)\frac{k!}{\ell!} x^{k-\ell}y^{\ell}z. \end{equation}
\Elm
\Bdm
Straightforward by induction on $k\in \mathbb N$. \Edm
\Bpo\label{po:gr}
Let us consider the filtration $F$ on the algebra $B$ given by the degree in $y$.
Then the differential of the Koszul complex $B\otimes_{B^e}K_w$ has degree $-1$ for
the total filtration degree. Moreover, the graded complex associated to the filtered
complex (\ref{shift}) is isomorphic to the Poisson homology complex of the Poisson
algebra $T$.
\Epo
\Bdm
As we have already seen, $gr_F(B)\simeq T$. As the monomials $x^iy^jz^k$, with
$i,j,k\in \mathbb N$, form a $\mathbb C$-basis of $B$ (Proposition
\ref{prop:basisofB}), the identifications of the spaces of Hochschild chains with
$B$ or $B^3$ and their analogs for the spaces of Poisson chains with $T$ and $T^3$
show that the spaces of Hochschild and Poisson chains are isomorphic. Now, we have to compare the images of $\tilde d_\ell$ and of $\delta_\ell$, $1\leq \ell\leq 3$. To do this, we consider the basis $(x^iy^jz^k)_{i,j,k}$ of $\mathbb C[x,y,z]$ and of $B$. In fact, using (\ref{eq:delta}), it is easy to compute the elements $\delta_\ell(x^iy^jz^k\, e)$, where $e$ is one of the symbols $d x$, $d y$ and $dz$ when $\ell=1$, $dy\wedge dz$, $dz\wedge dx$ and $dx\wedge dy$ when $\ell=2$, or $dx\wedge dy\wedge dz$, when $\ell=3$. Also, using lemma \ref{lm:form}, it is straightforward to  obtain $\tilde d_1(x^iy^jz^k\otimes x)$, $\tilde d_1(x^iy^jz^k\otimes y)$, $\tilde d_1(x^iy^jz^k\otimes z)$, $\tilde d_2(x^iy^jz^k\otimes r_1)$, $\tilde d_2(x^iy^jz^k\otimes r_2)$, $\tilde d_2(x^iy^jz^k\otimes r_3)$ and $\tilde d_3(x^iy^jz^k\otimes c(w))$ and to verify that
\begin{eqnarray*}
gr_F(\tilde d_1(x^i y^j z^k\otimes x)) &=&\delta_1(x^iy^jz^k dx)=jx^{i+2}y^{j-1}z^k,\\
gr_F(\tilde d_1(x^i y^j z^k\otimes y))&=&\delta_1(x^iy^jz^k dy)=(2k-i)x^{i+1}y^jz^k,\\
gr_F(\tilde d_1(x^i y^j z^k\otimes z))&=&\delta_1(x^iy^jz^k dz)=-2jx^{i+1}y^{j-1}z^{k+1},
\end{eqnarray*}
while:
\begin{eqnarray*}
gr_F(\tilde d_2(x^iy^jz^k \otimes r_1)) &\simeq& \delta_2(x^iy^jz^k dy\wedge dz)= \left( \begin{array}{c}2 x^iy^jz^{k+1}\\2jx^{i+1}y^{j-1}z^{k+1}
\\(2k-i+2)x^{i+1}y^jz^k \\ \end{array}\right),\\
gr_F(\tilde d_2(x^iy^jz^k \otimes r_2)) &\simeq& \delta_2(x^iy^jz^k dz\wedge dx) = \left( \begin{array}{c}-2j x^{i+1}y^{j-1}z^{k+1}\\0 \\-jx^{i+2}y^{j-1}z^k \\
\end{array}\right),\\
gr_F(\tilde d_2(x^iy^jz^k \otimes r_3)) &\simeq& \delta_2(x^iy^jz^k dx\wedge dy)= \left( \begin{array}{c}(i+2-2k) x^{i+1}y^{j}z^{k} \\jx^{i+2}y^{j-1}z^k \\0\\
\end{array}\right),
\end{eqnarray*}
and:
\begin{equation*}
gr_F(\tilde d_3(x^iy^jz^k\otimes c(w)))\simeq \delta_3(x^iy^jz^k dx\wedge dy\wedge
dz)
= \left(\begin{array}{c}jx^{i+2}y^{j-1}z^k\\ (2k-i)x^{i+1}y^j
z^k\\-2jx^{i+1}y^{j-1}z^{k+1} \end{array}\right).
\end{equation*}
This also shows that the differential $\tilde d_\ell$ has degree $-1$ for the total filtration degree. \ \Edm
%
%

%

In order to obtain the Poisson homology of the Poisson algebra $T$, we need to
compute the homology of a certain complex, called \emph{Koszul complex associated to
the polynomial} $\phi$. By definition, this complex is given by $0\rightarrow
\Omega^0(T)\rightarrow \Omega^1(T)\rightarrow \Omega^2(T)\rightarrow
\Omega^3(T)\rightarrow 0$, where the differential is the map $\wedge d\phi :
\Omega^k(T)\rightarrow \Omega^{k+1}(T)$. This complex will be denoted by $(\Omega^\bullet(T), \wedge d\phi)$ in the following, while for $0\leq p\leq 3$, the $p$-th homology space of
this complex will be denoted by $H_p^\phi(T)$.\\

\Rm Let us recall from \cite{pic:hom} that if $\varphi\in T$ is  a
weight-homogeneous polynomial with an isolated singularity at the origin, then the
homology of the Koszul complex associated to $\varphi$ is given by:
$H^\varphi_p(T)=\lbrace 0\rbrace$, for $p=0, 1, 2$, while $H_3^\varphi(T) =
\frac{T}{\lbrace\vec F\cdot \nabla\varphi\mid \vec F\in T^3\rbrace}$ is the
so-called Milnor algebra associated to $\varphi$ and is, in this case, a finite
dimensional vector space. In the following lemma, we see that the homology of the Koszul complex
associated to $\phi=-x^2z$ (admitting a non-isolated singularity at the origin) does
not satisfy the same properties.
\Blm\label{lm:koszul-phi}
Let $\phi = -x^2z$. The homology of the Koszul complex $(\Omega^\bullet(T), \wedge d\phi)$
associated to $\phi$ is given by:
\begin{eqnarray*} H^\phi_0(T) &=& \;\lbrace 0\rbrace ;\\
H^\phi_1(T) &=&  \mathbb C[y,z] \left(2z\, dx + x\, dz\right) ;\\ H^\phi_2(T) &=& \mathbb C[y,z] \left(x\, dz+2z\, dx\right)\wedge dy\oplus \left(x\mathbb C[y]\oplus \mathbb
C[y,z]\right) dz\wedge dx ;\\
H^\phi_3(T) &=& \left(x\mathbb C[y]\oplus \mathbb C[y,z]\right) dx\wedge dy\wedge dz.
\end{eqnarray*}
\Elm
\Bdm
The fact that $H^\phi_0(T)\simeq \lbrace 0\rbrace$ is clear. In order to compute $H^\phi_1(T)$, $H^\phi_2(T)$ and $H^\phi_3(T)$, we first point out that there exists an isomorphism of complexes between $(\Omega^\bullet(T), \wedge d\phi)$ and the tensor product of two analogous but simpler complexes:
\begin{equation}{\label{eq:isomcompl}}
(\Omega^\bullet(T), \wedge d\phi) \simeq (\Omega^\bullet(\mathbb C[y]), 0) \otimes (\Omega^\bullet(\mathbb C[x, z]), \wedge d\phi).
\end{equation}
Moreover, denoting $V=\mathbb C[x, z]$, the complex $(\Omega^\bullet(V), \wedge d\phi)$ can be written as follows:
\begin{equation*}
\setlength{\dgARROWLENGTH}{0.8cm}
\begin{diagram}
\node{\Omega^0(V)}\arrow{e,t}{\delta^0_V}\node{\Omega^1(V)}\arrow{e,t}{\delta^1_V}\node{\Omega^2(V)}\arrow{e,t}{}\node{0}
\end{diagram}
\end{equation*}
where $\Omega^0(V)=V$, $\Omega^1(V)=V\, dx\,\oplus V\, dz$ and $\Omega^2(V) = V\, dz\wedge dx$, while for every $(F,G) \in V^2$, $\delta^0_V(F)= F(2xz\, dx+x^2\, dz)$ and $\delta^1_V(F\, dx+G\, dz)=(Fx^2-2xzG)\, dz\wedge dx$.
Then, it is clear that the $0$-th cohomology space associated to $(\Omega^\bullet(V), \wedge d\phi)$ is $\lbrace 0\rbrace$ and the $2$-nd cohomology space is given by $\mathbb C[x,z]/ (x^2, 2xz) \simeq \left(\mathbb C[z]\oplus \mathbb C \cdot x\right)\, dz\wedge dx$.

In order to compute the first cohomology space associated to $(\Omega^\bullet(V), \wedge d\phi)$, let $F\, dx+G\, dz$  be a $1$-cochain of $(\Omega^\bullet(V), \wedge d\phi)$. This element is a cocycle if and only if $Fx^2-2xzG=0$, i.e., if and only if there exists $H\in V$ such that $F= 2zH$ and $G=xH$. Now, writing $H$ as $H=xK+L$, with $K\in V$ and $L\in \mathbb C[z]$, one has $F\, dx+G\, dz = \delta^0_V(K) + L\left(2z\, dx +x\, dz\right)$. Finally, for every $R\in \mathbb C[z]$, if $R\left(2z\, dx +x\, dz\right)\in \mathrm{Im}(\delta^0_V)$, then necessarily $R=0$. This allows us to conclude that the first cohomology space associated to the Koszul complex of $V$ can be written as $\mathbb C[z]\left(2z\, dx +x\, dz\right)$.

Since we are working over a field, from the isomorphism (\ref{eq:isomcompl}) and the K\"unneth formula, we immediately get $H^\phi_1(T)$, $H^\phi_2(T)$ and $H^\phi_3(T)$.
\Edm\\
%

\def\C{\mathbb{C}}

\Rm In \cite{pic:hom}, for every weight-homogeneous polynomial $\varphi\in \mathbb
C[x,y,z]$ admitting an isolated singularity, it is shown that $\lbrace \vec F\in T^3\mid
\nabla\varphi\cdot(\nabla\times\vec F)=0\rbrace = \lbrace \nabla G+ H\nabla
\varphi\mid G, H\in T\rbrace$, using the fact that the
Koszul complex associated to $\varphi$ is exact in degree $2$. In the following lemma, this is not true anymore if
$\varphi$ is replaced by $\phi = -x^2z$. \\

In the following, we will say that an element $\vec F=(F_1, F_2, F_3)\in T^3$ is
homogeneous of degree $d\in \mathbb Z$, if $F_1$, $F_2$, and $F_3$ are three
homogeneous polynomials of the same degree $d$.  (Notice that we use the convention that a polynomial of degree $d<0$ is zero).
\Blm\label{lm:eq-phi-rot}
Let $\phi=-x^2z\in T=\mathbb C[x,y,z]$. If $\vec F\in T^3$ (homogeneous or not) is such that  $\nabla\phi\cdot(\nabla\times\vec F)=0$, then there exist polynomials $G, H\in T$,  $C, P\in \C[y,z]$ and $f\in \mathbb C[\phi]$, satisfying  %
 \begin{equation}\label{eq:eqCP}
 \frac{\partial C}{\partial y} = P+2z\frac{\partial P}{\partial z}
 \end{equation}
 and
 \begin{equation}\label{eq:form2}
 \vec F= \nabla G+H\nabla \phi
+\left(\begin{smallmatrix}z\\0\\-x\end{smallmatrix}\right)C+\left(\begin{smallmatrix}2yz\\-3xz\\xy\end{smallmatrix}\right)P
+ f\left(\begin{smallmatrix}xz\\0\\-x^2\end{smallmatrix}\right).
 \end{equation}
 If $\vec F$ is homogeneous of degree $n\in \mathbb N$, then the polynomials $G, H, C, P$ and $f$ can be taken to be homogeneous of degrees $n+1$, $n-2$, $n-1$, $n-2$ and $n-2$, respectively. In particular, if $n\not\in 2+3\mathbb N$, then $f$ can be taken to be zero, and if $n=2+3k$, with $k\in \mathbb N$, then $f$ can be taken to be $f= \alpha (x^2z)^k$, with $\alpha\in \mathbb C$.
  Conversely, elements $\vec F\in T^3$ of the form (\ref{eq:form2}) satisfy the equation $\nabla\phi\cdot(\nabla\times \vec F)=0$.
\Elm
\Bdm
First of all, the last statement is straightforward.

Suppose that $\vec F$ is homogeneous of degree $n\in \mathbb N$ and satisfying $\nabla\phi\cdot(\nabla\times\vec F)=0$.
We will prove this result by recursion on $n\in \mathbb N$. As we have to distinguish whether $n\in 2+3\mathbb N$ or not, we first have to show
the desired result, for $n=0, 1, 2$.\\

$\circ$ If $n=0$, then there exist $a, b, c\in \C$, such that $\vec F=(a,b,c)$, and
it is clear that $\vec F=\nabla G$, with $G=ax+by+cz\in T$.\\

$\circ$ If $n=1$, there exist $a_i, b_i, c_i\in \C$, for $i=1, 2, 3$, such that $\vec F=
\left(\begin{smallmatrix}a_1x+b_1y+c_1z\\a_2x+b_2y+c_2z\\a_3x+b_3y+c_3z\end{smallmatrix}\right)$
and the condition $\nabla\phi\cdot(\nabla\times\vec F)=0$ is equivalent to:
$b_3=c_2$ and $a_2=b_1$. Using this, it is easy to verify that we can write $$
\vec F =  \nabla G+\left(\begin{smallmatrix}zC\\0\\-xC\end{smallmatrix}\right),
$$
with $G= \frac{1}{2}\left(a_1x^2+2b_1xy+(c_1+a_3)xz +b_2y^2+2c_2yz+c_3z^2\right)\in
T$ and $C=\frac{1}{2}(c_1-a_3)$.
Note also that $\frac{\partial C}{\partial y}=0$, so that the equation
(\ref{eq:eqCP}) is satisfied (here $P=0$).\\

$\circ$ If $n=2$, there exist $a_{ij}, b_{ij}, c_{ij}\in \C$, for $1\leq i,j\leq 3$,
such that $$
\vec F= \left(\begin{smallmatrix}a_{11}x^2+a_{12}xy+a_{13}xz +
a_{23}yz+a_{22}y^2+a_{33}z^2\\b_{11}x^2+b_{12}xy+b_{13}xz +
b_{23}yz+b_{22}y^2+b_{33}z^2\\c_{11}x^2+c_{12}xy+c_{13}xz +
c_{23}yz+c_{22}y^2+c_{33}z^2\end{smallmatrix}\right).
$$ The condition $\nabla\phi\cdot(\nabla\times\vec F)=0$ is equivalent to the following
identities: $2c_{12}-b_{13}-a_{23}=0$, $2c_{22}-b_{23}=0$, $c_{23}-2b_{33}=0$,
$2b_{11}-a_{12}=0$ and $b_{12}-2a_{22}=0$. Using this, it is straightforward to
verify that $$
\vec F = \nabla
G+\left(\begin{smallmatrix}zC+2yzP\\-3xzP\\-xC+xyP\end{smallmatrix}\right) +\alpha
\left(\begin{smallmatrix}xz\\0\\-x^2\end{smallmatrix}\right),
$$
with %
\begin{eqnarray*}
\lefteqn{G=
\frac{1}{3}a_{11}x^3+\frac{1}{2}a_{12}x^2y+\frac{1}{3}(a_{13}+c_{11})x^2z+\frac{1}{2}(a_{23}+b_{13})xyz + a_{22}xy^2}\\
&&\qquad \qquad \qquad +\frac{1}{3}(a_{33}+c_{13})xz^2+\frac{1}{2}b_{23}y^2z+b_{33}yz^2
+\frac{1}{3}b_{22}y^3+\frac{1}{3}c_{33}z^3,
\end{eqnarray*}
 $\alpha=\frac{1}{3}(a_{13}-2c_{11})$,
$C=\frac{1}{6}(a_{23}-b_{13})y+\frac{1}{3}(2a_{33}-c_{13})z$ and
$P=\frac{1}{6}(a_{23}-b_{13})$. Notice that the equation (\ref{eq:eqCP}) is clearly
satisfied in this case.\\
 $\circ$ Let now $m\in \mathbb N$, such that $m\geq 3$, and suppose the lemma is
proved for all $n\in \mathbb N$ such that $n<m$. Let now $\vec F\in T^3$ be a
homogeneous element of degree $m$, satisfying $\nabla\phi\cdot(\nabla\times\vec
F)=0$.

According to Lemma \ref{lm:koszul-phi}, this hypothesis implies that there exist
homogeneous elements: $\vec K\in T^3$,  $P, C\in \C[y,z]$ and $E\in \C[y]$ such that
$$
\nabla\times \vec F = \vec K\times \nabla \phi +
\left(\begin{smallmatrix}x\\0\\-2z\end{smallmatrix}\right)P +
\left(\begin{smallmatrix}0\\1\\0\end{smallmatrix}\right)(xE+C),
$$
with the degrees of $\vec K$, $P$, $E$ and $C$ respectively equal to $m-3$, $m-2$,
$m-2$ and $m-1$.
Computing the divergence of this, we obtain
$$
0= \mathrm{Div}(\nabla\times \vec F) = (\nabla\times \vec K)\cdot \nabla \phi -P
-2z\frac{\partial P}{\partial z}+xE'+\frac{\partial C}{\partial y},
$$
where $E'=\frac{d E}{dy}$ and where we have used that $\mathrm{Div}(\vec K\times
\nabla \phi) = (\nabla\times \vec K)\cdot \nabla \phi$. This implies that $x$
divides the polynomial $-P -2z\frac{\partial P}{\partial z}+\frac{\partial
C}{\partial y}$ which lies in $\C[y,z]$, so that it is zero:
$$
\frac{\partial C}{\partial y} = P +2z\frac{\partial P}{\partial z}.
$$
It remains that $ (\nabla\times \vec K)\cdot \nabla \phi+xE'=0$, which gives $E'\in
(xT+zT)\cap \C[y]=\lbrace 0\rbrace$, and $E'=0$. This fact, together with the
hypothesis that $E$ is supposed to be homogeneous of degree $m-2$, and $m\geq 3$,
imply $E=0$. We have obtained:
$$
\nabla\times \vec F = \vec K\times \nabla \phi +
\left(\begin{smallmatrix}xP\\C\\-2zP\end{smallmatrix}\right).
$$
Now, $\vec K$ is of degree $m-3<m$ and satisfies $(\nabla\times \vec K)\cdot \nabla
\phi = 0$, so that we can apply the recursion hypothesis to obtain the existence of
homogeneous elements $G, H\in T$, $D, Q\in \C[y,z]$, $\alpha\in \C$, and $k\in
\mathbb N$, such that %
\begin{equation}\label{eq:DQyz}
\frac{\partial D}{\partial y}=Q+2z\frac{\partial Q}{\partial z},
\end{equation}
and %
\begin{equation*}
\vec K = \nabla G+H\nabla \phi
+\left(\begin{smallmatrix}zD+2yzQ\\-3xzQ\\-xD+xyQ\end{smallmatrix}\right) +\alpha
(x^2z)^k\left(\begin{smallmatrix}xz\\0\\-x^2\end{smallmatrix}\right).
\end{equation*}
The polynomials $G$, $H$, $D$ and $Q$ are respectively of degree $m-2$, $m-5$, $m-4$
and $m-5$.
Notice that, by hypothesis, $\alpha$ is supposed to be zero, except if $m-3 = 2+3k$.
We now compute
\begin{eqnarray*}
\vec K\times \nabla\phi  &=& \nabla G\times\nabla\phi  + \left(\begin{smallmatrix}3x^3zQ\\3x^2zD\\-6x^2z^2Q\end{smallmatrix}\right) +\alpha
(x^2z)^k\left(\begin{smallmatrix}0\\3x^3z\\0\end{smallmatrix}\right).
\end{eqnarray*}
Denoting by $\vec e$ the so-called Euler vector $\vec e:=(x,y,z)\in T^3$, we
use the following general result (Proposition 3.5 in \cite{pic:hom}), which is due
to the exactness of the De Rham complex of $\C[x,y,z]$: if $\vec A=(A_1, A_2,
A_3)\in T^3$ is a  homogeneous element of degree $d\in \mathbb N$, such that
$\mathrm{Div}(\vec A)=0$, then the Euler formula $\nabla A_i\cdot \vec e= d\,A_i$
($1\leq i\leq 3$), implies that $(d+2) \vec A =\nabla\times(\vec A\times \vec e)$.

As %
\begin{eqnarray*}
\mathrm{Div}\left(\left(\begin{smallmatrix}3x^3zQ\\3x^2zD\\-6x^2z^2Q\end{smallmatrix}\right)
+\alpha (x^2z)^k\left(\begin{smallmatrix}0\\3x^3z\\0\end{smallmatrix}\right)\right)
&=& \mathrm{Div}\left(\vec K\times \nabla\phi  -\nabla G\times\nabla\phi  \right)\\
&=&\left(\nabla\times\left(\vec K-\nabla G\right)\right)\cdot \nabla \phi\\
&=&0,
\end{eqnarray*}
we can write $\vec K\times \nabla\phi= \nabla\times \vec L$, where $\vec L$ is given
by:
\begin{eqnarray*}
\vec L
 &=& G\nabla\phi  +
\frac{1}{(m+1)}\left(\begin{smallmatrix}3x^3zQ\\3x^2zD\\-6x^2z^2Q\end{smallmatrix}\right)\times\vec
e +\frac{1}{m+1}\alpha
(x^2z)^k\left(\begin{smallmatrix}0\\3x^3z\\0\end{smallmatrix}\right)\times\vec e\\
&=& G\nabla\phi  + \frac{1}{(m+1)}\left(\begin{smallmatrix}3x^2z^2D
+6x^2yz^2Q\\-9x^3z^2Q\\3x^3yzQ-3x^3zD\end{smallmatrix}\right) +\frac{3}{m+1}\alpha
(x^2z)^{k+1}\left(\begin{smallmatrix}xz\\0\\-x^2\end{smallmatrix}\right).\end{eqnarray*}
Moreover, we have also %
\begin{equation*}
\left(\begin{smallmatrix}xP\\C\\-2zP\end{smallmatrix}\right) =\frac{1}{m+1}
\nabla\times\left(\left(\begin{smallmatrix}xP\\C\\-2zP\end{smallmatrix}\right)\times\vec
e\right)
= \frac{1}{m+1}
\nabla\times\left(\left(\begin{smallmatrix}zC+2yzP\\-3xzP\\-xC+xyP\end{smallmatrix}\right)\right).
\end{equation*}
We finally obtain %
\begin{equation*}
\nabla\times \vec F =  \nabla\times \left(\vec L+ \frac{1}{m+1}
\left(\begin{smallmatrix}zC+2yzP\\-3xzP\\-xC+xyP\end{smallmatrix}\right)\right).
\end{equation*}
This allows us to apply another general result (Proposition 3.5 in \cite{pic:hom}):
if $\vec A=(A_1, A_2, A_3)\in T^3$ is a  homogeneous element of degree $d\in \mathbb
N$, such that $\nabla\times \vec A=0$, then the Euler formula $\nabla A_i\cdot \vec
e= d\,A_i$ ($1\leq i\leq 3$), implies that $(d+1) \vec A =\nabla(\vec A\cdot \vec
e)$.

This implies that there exists a homogeneous element $S\in T$ (of degree $m+1$),
such that
\begin{eqnarray*}
\vec F &=& \nabla S +\vec L+ \frac{1}{m+1}
\left(\begin{smallmatrix}zC+2yzP\\-3xzP\\-xC+xyP\end{smallmatrix}\right)\\
& =& \nabla S +G\nabla\phi  + \frac{1}{(m+1)}\left(\begin{smallmatrix}3x^2z^2D
+6x^2yz^2Q\\-9x^3z^2Q\\3x^3yzQ-3x^3zD\end{smallmatrix}\right) \\
&&\\
&&\qquad\qquad\qquad+\frac{3}{m+1}\alpha
(x^2z)^{k+1}\left(\begin{smallmatrix}xz\\0\\-x^2\end{smallmatrix}\right)+\frac{1}{m+1}
\left(\begin{smallmatrix}zC+2yzP\\-3xzP\\-xC+xyP\end{smallmatrix}\right) .
\end{eqnarray*}
Let now $V:=\frac{3}{(2m-7)}xz(D+2yQ)$ (as $2m-7\not=0$). The polynomial $V$ is
homogeneous of degree $m-2$.
Now, using (\ref{eq:DQyz}) and the Euler formulas $y\frac{\partial D}{\partial
y}+z\frac{\partial D}{\partial z}=(m-4)D$ and $y\frac{\partial Q}{\partial
y}+z\frac{\partial Q}{\partial z}=(m-5)Q$, it is straightforward to verify that we have :
\begin{eqnarray*}
\nabla\left(-\frac{3}{m+1} x^2z V\right)-V\nabla\phi = \frac{1}{(m+1)}\left(\begin{smallmatrix}3x^2z^2D
+6x^2yz^2Q\\-9x^3z^2Q\\3x^3yzQ-3x^3zD\end{smallmatrix}\right).
\end{eqnarray*}
Denoting by $\tilde\alpha:=\frac{3}{m+1}\alpha$, $\tilde G := G-V$ and by $\tilde
S:=S-\frac{3}{(m+1)}x^2zV$, $\tilde C=\frac{1}{m+1}C$ and $\tilde P=\frac{1}{m+1}P$,
we can write:
\begin{eqnarray*}
\vec F &=&  \nabla\tilde S +\tilde G\nabla\phi  +\tilde \alpha
(x^2z)^{k+1}\left(\begin{smallmatrix}xz\\0\\-x^2\end{smallmatrix}\right)+\left(\begin{smallmatrix}z\tilde C+2yz\tilde P\\-3xz\tilde P\\-x\tilde C+xy\tilde
P\end{smallmatrix}\right) .
\end{eqnarray*}
Moreover, $\tilde S$ and $\tilde G$ are homogeneous polynomials and we have already
seen that $C$ and $P$ satisfy the identity (\ref{eq:eqCP}) so that this identity is
also satisfied by $\tilde C$ and $\tilde P$.
\Edm\\
In the sequel, we will several times need the following technical result.%
\Blm\label{lm:phidivides}
If $r\in \mathbb N$, $c\in \C$ and if $K\in \C[x,y,z]$ is a polynomial such that $c x\phi^r
\left(\begin{smallmatrix}0\\1\\0\end{smallmatrix}\right) = \nabla
K\times\nabla\phi$, then $K\in \C[\phi]$ and $c$ is necessarily $0$.
\Elm
\Bdm
The hypothesis on $K$ is equivalent to $\frac{\partial K}{\partial y}=0$ and
$c\phi^r = -2z\frac{\partial K}{\partial z}+x\frac{\partial K}{\partial x}$, hence
$K\in \C[x,z]$.
Denote by $\mathcal D$ the operator $\mathcal D = -2z\frac{\partial}{\partial z} + x\frac{\partial}{\partial x}$. For every $i,j\in \mathbb N$, $\mathcal D(x^iz^j)=(i-2j) x^iz^j$, so that $\mathcal D$ is diagonalizable on $\mathbb C[x,z]$ and its kernel is $\mathbb C[x^2z]$. As $\mathcal D(K) = c\phi^r$, we have $$
0 = \mathcal D(c\phi^r) = \mathcal D^2 (K),
$$
hence $K\in \ker (\mathcal D) = \mathbb C[\phi]$. Finally, the fact that $K$ is in $\mathbb C[\phi]$ implies $\nabla
K\times\nabla\phi = 0$, and $c=0$.
 \Edm\\

\Bcr\label{cr:equation}
Let $\phi=-x^2z\in T=\mathbb C[x,y,z]$. We have %
\begin{eqnarray*}
\frac{\lbrace \vec F\in T^3\mid \nabla\phi\cdot(\nabla\times \vec
F)=0\rbrace}{\lbrace \nabla G+ H\nabla\phi\mid G, H\in T\rbrace} &\simeq&
 \C[\phi] \left(\begin{smallmatrix}xz\\0\\-x^2\end{smallmatrix}\right)\oplus   \C[z] \left(\begin{smallmatrix}z^2\\0\\-xz\end{smallmatrix}\right)\\
  &&\oplus \renewcommand{\arraystretch}{0.7}\bigoplus_{\begin{array}{c}\scriptstyle
n\in \mathbb N\\\scriptstyle 1\leq k\leq n+1\end{array}}
\C\,\left(\begin{smallmatrix}(2n+3)\,yz\\-3k\,xz\\(-2n+3(k-1))\,xy\end{smallmatrix}\right)y^{k-1}z^{n+1-k}.
\end{eqnarray*}
\Ecr
\Bdm
Lemma \ref{lm:eq-phi-rot} gives:
\begin{eqnarray*}
\lefteqn{\frac{\lbrace \vec F\in T^3\mid \nabla\phi\cdot(\nabla\times \vec
F)=0\rbrace}{\lbrace \nabla  G+H\nabla\phi\mid G, H\in T\rbrace}}\\
&&\\
&=& \left
\lbrace\left(\begin{smallmatrix}zC+2yzP\\-3xzP\\-xC+xyP\end{smallmatrix}\right)
+\alpha (x^2z)^k\left(\begin{smallmatrix}xz\\0\\-x^2\end{smallmatrix}\right)\mid
\alpha\in \C, k\in \mathbb N, \right.\\
&&\qquad \qquad \qquad\left. C, P\in \C[y,z]\hbox{ satisfying }
\scriptstyle\frac{\partial C}{\partial y} = P+2z\frac{\partial P}{\partial
z}\right\rbrace.
\end{eqnarray*}
Fix now $n\in \mathbb N$. Let $C, P\in \C[y,z]$ be homogeneous polynomials
satisfying $\frac{\partial C}{\partial y} = P+2z\frac{\partial P}{\partial z}$. We
suppose that $P$ is of degree $n$, so that $C$ is zero or of degree equal to $n+1$.
We write $C$ and $P$ as
$$
P = \sum_{k=0}^n a_k y^k z^{n-k}, \hbox{ and } C= \sum_{k=0}^{n+1} b_k y^k z^{n+1-k},
$$
where $a_k, b_k\in \C$. Now, compute
\begin{equation*}
0=\frac{\partial C}{\partial y} -P-2z\frac{\partial P}{\partial z} = \sum_{k=0}^n
\left((k+1) b_{k+1}- (2(n-k)+1) a_k\right) y^kz^{n-k},
\end{equation*}
so that, necessarily, for all $k=1, \dots, n+1$, we have $b_k = \frac{(2(n-k)+3)}{k}
a_{k-1}$. We then can write
\begin{eqnarray*}
\left(\begin{smallmatrix}zC+2yzP\\-3xzP\\-xC+xyP\end{smallmatrix}\right)  &=&
 b_0\left(\begin{smallmatrix}z^{2}\\0\\-xz\end{smallmatrix}\right)z^n  +\sum_{k=1}^{n+1}  \frac{1}{k}a_{k-1} \left(\begin{smallmatrix}(2n+3)yz\\-3k
xz\\(-2n+3(k-1))xy\end{smallmatrix}\right)y^{k-1} z^{n-k+1}\\
 &\in &\C[z] \left(\begin{smallmatrix}z^2\\0\\-xz\end{smallmatrix}\right)\,  +\, \renewcommand{\arraystretch}{0.7}\sum_{k=1}^{n+1}
\,\C\,\left(\begin{smallmatrix}(2n+3)\,yz\\-3k\,xz\\(-2n+3(k-1))\,xy\end{smallmatrix}\right)y^{k-1}z^{n+1-k}
\end{eqnarray*}
We finally have shown that
\begin{eqnarray*}
\frac{\lbrace \vec F\in T^3\mid \nabla\phi\cdot(\nabla\times \vec
F)=0\rbrace}{\lbrace \nabla G+H\nabla\phi\mid G, H\in T\rbrace} &\simeq&
 \C[\phi] \left(\begin{smallmatrix}xz\\0\\-x^2\end{smallmatrix}\right)+
  \C[z] \left(\begin{smallmatrix}z^2\\0\\-xz\end{smallmatrix}\right)\\
  &+&\!\!\!\! \renewcommand{\arraystretch}{0.7}\sum_{\begin{array}{c}\scriptstyle n\in
\mathbb N\\\scriptstyle 1\leq k\leq n+1\end{array}}
\C\,\left(\begin{smallmatrix}(2n+3)\,yz\\-3k\,xz\\(-2n+3(k-1))\,xy\end{smallmatrix}\right)y^{k-1}z^{n+1-k},
\end{eqnarray*}
and it remains to show that this sum is a direct one. To do this, it suffices to
show that each homogeneous component of this sum is a direct sum. Notice that an
element of the space of the right hand side of the previous equation is at least of
degree $2$, so that, we fix $n\in \mathbb N$, and we consider the degree $n+2$
component of the previous sum. Let $\alpha, c\in \C$ and $r\in \mathbb N$, and for
all $1\leq k\leq n+1$, we consider $a_k\in \C$. Suppose that there exist
homogeneous elements $G, H\in T$, of respective degrees equal to $n+3$ and $n$,
satisfying:
\begin{equation}\label{eq:direct}
\alpha\phi^r\left(\begin{smallmatrix}xz\\0\\-x^2\end{smallmatrix}\right) + cz^n \left(\begin{smallmatrix}z^2\\0\\-xz\end{smallmatrix}\right)
+\sum_{k=1}^{n+1} a_k \left(\begin{smallmatrix}(2n+3) y^k
z^{n+2-k}\\-3kxy^{k-1}z^{n+2-k}\\(-2n+3(k-1))\,xy^kz^{n+1-k}\end{smallmatrix}\right)
= \nabla G+H\nabla \phi.
\end{equation}
Notice that $\alpha=0$ if $n\not= 3r$.
Applying the curl operator to this identity allows us to obtain
$$
-\sum_{k=1}^{n+1} k(2n+6)\, a_k\, y^{k-1} z^{n+2-k} = 2xz\frac{\partial H}{\partial
y}, $$
which implies that $a_k=0$, for  all $1\leq k \leq n+1$. This, together with
(\ref{eq:direct}) imply that $\frac{\partial G}{\partial y}=0$, so that $G\in
\C[x,z]$, while, this together with the result obtained by applying the curl
operator to the previous identity (\ref{eq:direct}), give:
$$
3(r+1)\alpha \,x\phi^r + (n+3) cz^{n+1} = -2xz\frac{\partial H}{\partial z}
+x^2\frac{\partial H}{\partial x}.
$$
This shows that $x$ divides $(n+3)c z^{n+1}$, which means that $c=0$ and the
equation (\ref{eq:direct}) becomes: %
\begin{equation}\label{eq:eqalpha}
\alpha\phi^r\left(\begin{smallmatrix}xz\\0\\-x^2\end{smallmatrix}\right)  = \nabla G+H\nabla \phi.
 \end{equation}
 Applying to this equation the cross product with $\nabla \phi$ leads to: $-3\alpha x\phi^{r+1}
\left(\begin{smallmatrix}0\\1\\0\end{smallmatrix}\right) = \nabla
G\times\nabla\phi$. This, together with lemma \ref{lm:phidivides}, give $\alpha=0$, which is what remains to prove.
\Edm\\

We now determine the Poisson homology of the Poisson algebra $(T, \lbrace\cdot,
\cdot\rbrace)$. %
\Bpo\label{po:poissonhomol}
 Set $T=\mathbb{C}[x,y,z]$ and consider $\phi=-x^2z\in T$.  The algebra $T$ becomes
a Poisson algebra when equipped with the Poisson bracket $ \lbrace\cdot,
\cdot\rbrace$, defined by:
$$
\lbrace y,z\rbrace = \frac{\partial \phi}{\partial x}=-2xz, \quad \lbrace z,x\rbrace
= \frac{\partial \phi}{\partial y}=0,\quad \lbrace x,y\rbrace = \frac{\partial
\phi}{\partial z}=-x^2.  $$

Using the identifications $\Omega^1(T)\simeq T^3$, $\Omega^2(T)\simeq T^3$,
$\Omega^3(T)\simeq T$ explained above, the Poisson homology spaces of the Poisson
algebra $T$ are given by:
\begin{eqnarray*}
HP_0(T) &\simeq& x\mathbb{C}[y]\oplus \mathbb{C}[y,z];\\
 {}\\
HP_1(T) &\simeq& \mathbb{C}[\phi]
\left(\begin{smallmatrix}xz\\0\\-x^2\end{smallmatrix}\right)\oplus\mathbb{C}[z]
\left(\begin{smallmatrix}z^2\\0\\-xz\end{smallmatrix}\right)\oplus  \renewcommand{\arraystretch}{0.7}
  \bigoplus_{\begin{array}{c}\scriptstyle n\in \mathbb N^*\\\scriptstyle 0\leq k\leq
n\end{array}}
 \mathbb{C}  \left(\begin{smallmatrix}0\\ky^{k-1}z^{n-k}\\(n-k)y^kz^{n-1-k}\end{smallmatrix}\right)\\
&\oplus &\bigoplus_{n\in \mathbb
N}\,\mathbb{C}\left(\begin{smallmatrix}y^n\\n\,xy^{n-1}\\0\end{smallmatrix}\right)
\oplus  \renewcommand{\arraystretch}{0.7}
  \bigoplus_{\begin{array}{c}\scriptstyle n\in \mathbb N\\ \scriptstyle 1\leq k\leq
n+1\end{array}}\mathbb C\left(\begin{smallmatrix}(2n+3)\, yz\\-3k\,xz\\
(-2n+3(k-1))\,xy\end{smallmatrix}\right)y^{k-1}z^{n+1-k};
\end{eqnarray*}
\begin{eqnarray*}
HP_2(T) &\simeq& \mathbb C[\phi]
\left(\begin{smallmatrix}x\\y\\z\end{smallmatrix}\right)\,\oplus\,   \left(x\mathbb C[\phi]\oplus z\mathbb
C[z]\right)\left(\begin{smallmatrix}0\\1\\0\end{smallmatrix}\right)\\
 &&\qquad \qquad \qquad \qquad \qquad \qquad \oplus \renewcommand{\arraystretch}{0.7}
  \bigoplus_{\begin{array}{c}\scriptstyle n\in \mathbb N\\\scriptstyle 0\leq k\leq
n\end{array}}\mathbb C\left(\begin{smallmatrix}(k+1) x\\ (2(n-k)+1)\,y\\
-2(k+1)\,z\end{smallmatrix}\right)y^{k}z^{n-k};\\
  {}\\
  HP_3(T) &\simeq& \mathbb C[\phi].
\end{eqnarray*}
\Epo
\Bdm Remark that it is not possible to apply the results of Monnier for the computation of Poisson homology as we have seen for the third type (proof of Proposition \ref{quantum}). Actually, it is easy to check that the Poisson bracket derived from the potential $x^2z$ is not diagonalizable. To determine the Poisson homology spaces, notice that the polynomial $\phi$ is homogeneous of
degree $3$, so that, considering $T$ graded by the total degree of the polynomials,
the operator $\delta_p$ ($1\leq p\leq 3$) is homogeneous of degree $1$. This allows
one to determine the Poisson homology spaces, degree by degree.\\

\noindent\textit{The $0$-th Poisson homology space $HP_0(T)$.}\\
According to (\ref{eq:delta}), $$
HP_0(T)\simeq \frac{T}{\left\lbrace \delta_1(\vec G)\mid\vec G = (G_1, G_2,
G_3)\in T^3\right\rbrace},
$$
where $ \delta_1(\vec G)=2xz\left(\frac{\partial
G_2}{\partial z}-\frac{\partial G_3}{\partial y}\right)+x^2\left(\frac{\partial
G_1}{\partial y}-\frac{\partial G_2}{\partial x}\right)$.
Now let $G\in T$. We can write $G$ as $G=xF+H$ with $F\in T$ and $H\in \mathbb C[y,z]$. Then, $F$ can also be written as $F= xA+zB+C$, with $A, B\in T$ and $C\in \mathbb C[y]$, so that $$
G= x^2A+xzB+xC+H.
$$
Now, there exist $G_1\in T$ and $G_3\in T$, such that $A = \frac{\partial
G_1}{\partial y}$ and $B = -2\frac{\partial G_3}{\partial y}$, so that $$
G = x^2A + xzB +xC +H =-2xz\frac{\partial G_3}{\partial y} +x^2\frac{\partial
G_1}{\partial y} + xC+H
=\delta_1(\vec G)+xC+H,
$$
with $\vec G=(G_1, 0, G_3)\in T^3$.
This shows that $HP_0(T) \simeq x\mathbb C[y]+\mathbb C[y,z]$.
It is also clear
that this sum is a direct one, because an equality of the form $xC+H= \vec
F\cdot\nabla\phi$, with  $C\in \mathbb C[y]$, $H\in \mathbb C[y,z]$ and $\vec F\in
T^3$ implies that $x$ divides $H$, so $H=0$ and then it remains that $C \in
(xT+zT)\cap \mathbb C[y]$ which means that $C=0$.
 We
finally have obtained $HP_0(T)\simeq x\mathbb C[y]\oplus\mathbb C[y,z]$.\\

\noindent\textit{The third Poisson homology space $HP_3(T)$.}\\
First, it is easy to see that $\mathbb C[\phi]\simeq \mathbb C[\phi]dx\wedge dy
\wedge dz\subseteq HP_3(T)$.
Conversely, we consider an element $F\in HP_3(T)$, i.e., $F\in T\simeq
\Omega^3(T)$ satisfying $\nabla F\times \nabla(x^2z)=0$. Lemma \ref{lm:phidivides} (with $c=0$)
then implies that $F\in \mathbb C[\phi]$, so that $HP_3(T)\simeq \mathbb C[x^2z ] =\mathbb C[\phi]$.\\

\noindent\textit{The first Poisson homology space $HP_1(T)$.}\\
By (\ref{eq:delta}),
$$
HP_1(T) = \frac{\left\lbrace \vec F\in T^3\mid \nabla\phi\cdot(\nabla\times\vec F)=0
\right\rbrace}{\left\lbrace -\nabla\left(\vec G\cdot
\nabla\phi\right)+\mathrm{Div}(\vec G)\nabla \phi\mid \vec G\in T^3\right\rbrace}.
$$
Let $\vec F\in T^3$ be an element satisfying  $\nabla\phi\cdot(\nabla\times\vec
F)=0$. According to corollary \ref{cr:equation},  there exist $G,H\in T$ such that:
\begin{eqnarray*}
\vec F \in \nabla G+H\nabla \phi+ \C[\phi] \left(\begin{smallmatrix}xz\\0\\-x^2\end{smallmatrix}\right) + \C[z]
\left(\begin{smallmatrix}z^2\\0\\-xz\end{smallmatrix}\right)+ \sum_{n\in\mathbb
N}\sum_{k=1}^{n+1} \C\; \vec u_{n,k},
\end{eqnarray*}
where $\vec u_{n,k} =
\left(\begin{smallmatrix}(2n+3)\,yz\\-3k\,xz\\(-2n+3(k-1))\,xy\end{smallmatrix}\right)y^{k-1}z^{n+1-k}$, for $n\in \mathbb N$ and $1\leq k \leq n+1$.

Now, according to the determination of $HP_0(T)$, there exist $\vec L, \vec K\in
T^3$, $A, \tilde A\in \C[y]$ and $B, \tilde B\in \C[y,z]$ such that $$
G = \nabla \phi\cdot(\nabla\times \vec L) + xA +B, \quad \hbox{ and }\quad H =
\nabla \phi\cdot(\nabla\times \vec K) + x\tilde A +\tilde B.
$$
As $\delta_2(-\nabla\times\vec L) = \nabla\left( \nabla\phi\cdot(\nabla\times\vec
L)\right)$, and $\delta_2(\vec K\times \nabla\phi)= \left(\nabla\phi\cdot(\nabla\times \vec
K)\right)\nabla\phi$, we obtain
\begin{eqnarray*}
\vec F &\in& \delta_2(-\nabla\times\vec L+\vec K\times \nabla\phi)+\nabla (xA
+B)+(x\tilde A +\tilde B)\nabla\phi\\
&+& \C[\phi]  \left(\begin{smallmatrix}xz\\0\\-x^2\end{smallmatrix}\right) + \C[z]
\left(\begin{smallmatrix}z^2\\0\\-xz\end{smallmatrix}\right)+ \sum_{n\in\mathbb
N}\sum_{k=1}^{n+1} \C\; \vec u_{n,k}.
\end{eqnarray*}
Now, let us consider $A_1\in \C[y]$ satisfying $\frac{\partial A_1}{\partial
y}=\tilde A$ and $B_1\in \C[y,z]$ satisfying $\frac{\partial B_1}{\partial y}=\tilde
B$. It is then straightforward to verify that $$
\delta_2\left(\left(\begin{smallmatrix}-x\tilde B\\xA_1-2z\frac{\partial
B_1}{\partial z}\\2z\tilde B\end{smallmatrix}\right)\right) = (x\tilde A+\tilde
B)\nabla\phi,
$$
so that $$
\vec F \in \mathrm{Im}(\delta_2)+\nabla (xA +B) + \C[\phi] \left(\begin{smallmatrix}xz\\0\\-x^2\end{smallmatrix}\right) + \C[z]
\left(\begin{smallmatrix}z^2\\0\\-xz\end{smallmatrix}\right)+ \sum_{n\in\mathbb
N}\sum_{k=1}^{n+1} \C\; \vec u_{n,k}.
$$
This allows us to write %
\begin{eqnarray}\label{eq:HP1direct}
HP_1(T) &\simeq&
\C[\phi]\left(\begin{smallmatrix}xz\\0\\-x^2\end{smallmatrix}\right) +
\C[z]\left(\begin{smallmatrix}z^2\\0\\-xz\end{smallmatrix}\right) +
\sum_{n\in\mathbb N}\sum_{k=1}^{n+1} \C\; \vec u_{n,k} \nonumber\\
&+& \left\lbrace \nabla (xA +B) \mid A\in \C[y], B\in y\C[y,z]+z\C[y,z]\right\rbrace.
\end{eqnarray}
Let us show that this sum is a direct one. Suppose that there exist $P\in \C[X]$ a polynomial in one variable, $Q\in \C[z]$,
$A\in \C[y]$, and a polynomial $B\in y\C[y,z]+z\C[y,z]$ (i.e., $B\in \C[y,z]$ with no constant term), $\vec K\in T^3$, and for
every $n\in \mathbb N$ and every $1\leq k\leq n+1$, suppose that $\alpha^n_k\in \C$
are constants, such that:
\begin{eqnarray*}
\lefteqn{ \nabla (xA+B)+ P(\phi) \left(\begin{smallmatrix}xz\\0\\-x^2\end{smallmatrix}\right) + Q(z)
\left(\begin{smallmatrix}z^2\\0\\-xz\end{smallmatrix}\right)+ \sum_{n\in \mathbb
N}\sum_{k=1}^{n+1} \alpha^n_k\, \vec u_{n,k}}\\
&&\qquad \qquad \qquad \qquad \qquad\qquad = - \nabla(\vec
K\cdot\nabla\phi)+\mathrm{Div}(\vec K)\nabla\phi.
\end{eqnarray*}
This implies that $$
 P(\phi)  \left(\begin{smallmatrix}xz\\0\\-x^2\end{smallmatrix}\right) + Q(z)
\left(\begin{smallmatrix}z^2\\0\\-xz\end{smallmatrix}\right)+ \sum_{n\in\mathbb
N}\sum_{k=1}^{n+1} \alpha^n_k \vec u_{n,k} \in \lbrace K\nabla\phi+ \nabla L\mid
K,L\in T\rbrace.
$$
According to Corollary \ref{cr:equation}, necessarily $P=0$, $Q=0$ and $\alpha^n_k=0$, for all $n\in \mathbb N$ and all $1\leq
k\leq n+1$. It then remains $\nabla (xA+B+\vec K\cdot\nabla\phi) = \mathrm{Div}(\vec
K)\nabla\phi$, which implies that $\nabla (xA+B+\vec K\cdot\nabla\phi)
\times\nabla\phi=0$ and according to the determination of the space $HP_3(T)$, this
gives the existence of polynomial in one variable $R\in \C[X]$, such that $xA+B+\vec K\cdot\nabla\phi = R(\phi)$. As there is no constant term in $B$ and for
degree reason, necessarily $\phi=-x^2z$ divides $R(\phi)$ and $x$ divides $B\in
\C[y,z]$, which implies that $B=0$. Moreover, last equation allows us to obtain $
A\in \C[y]\cap (xT+zT)=\lbrace 0\rbrace$.
This allows us to conclude that the sum in (\ref{eq:HP1direct}) is a direct one,
i.e., we can write
\begin{eqnarray*}
HP_1(T) &\simeq&
\C[\phi]\left(\begin{smallmatrix}xz\\0\\-x^2\end{smallmatrix}\right) \oplus
\C[z]\left(\begin{smallmatrix}z^2\\0\\-xz\end{smallmatrix}\right)
\oplus\bigoplus_{n\in\mathbb N}\bigoplus_{k=1}^{n+1} \C\; \vec u_{n,k}\\
&\oplus& \left\lbrace \nabla (xA +B) \mid A\in \C[y], B\in
y\C[y,z]+z\C[y,z]\right\rbrace.
\end{eqnarray*}
Finally, it is clear that %
\begin{eqnarray*}
 \lefteqn{\left\lbrace \nabla (xA +B) \mid A\in \C[y], B\in
y\C[y,z]+z\C[y,z]\right\rbrace = }\\
 &&\qquad\qquad\bigoplus_{n\in \mathbb
N}\,\mathbb{C}\left(\begin{smallmatrix}y^n\\n\,xy^{n-1}\\0\end{smallmatrix}\right)\oplus
\bigoplus_{n\in \mathbb N^*}\bigoplus_{0\leq k\leq n}
 \mathbb{C}  \left(\begin{smallmatrix}0\\ky^{k-1}z^{n-k}\\(n-k)y^kz^{n-1-k}\end{smallmatrix}\right),
\end{eqnarray*}
and this finishes the determination of $HP_1(T)$.\\
\noindent\textit{The second Poisson homology space $HP_2(T)$.}\\
According to (\ref{eq:delta}), we have
$$
HP_2(T) = \frac{\left\lbrace\vec F\in T^3\mid -\nabla(\vec F\cdot \nabla
\phi)+\mathrm{Div}(\vec F)\nabla\phi=0 \right\rbrace}{\left\lbrace \nabla G\times
\nabla\phi\mid G\in T\right\rbrace}.
$$
Let $\vec F\in T^3$ be a homogeneous element of degree $n\in \mathbb N$, satisfying %
\begin{equation}\label{eq:HP2}
-\nabla(\vec F\cdot \nabla \phi)+\mathrm{Div}(\vec F)\nabla\phi=0.
\end{equation}
This implies that $\nabla(\vec F\cdot \nabla \phi)\times \nabla\phi=0.$ This, together with the determination of $HP_3(T)$ gives the existence of
$\alpha\in \C$ and $r\in \mathbb N$ such that $\vec F\cdot \nabla \phi = \alpha
\phi^r$. Notice that if $r=0$, then for degree reasons, $\alpha=0$.
According to the Euler formula $\nablaæ\phi\cdot \vec e = 3\phi$ (where we
recall that $\vec e = (x,y,z)\in T^3$), we get $\vec F\cdot \nabla \phi = \frac{
\alpha}{3}\phi^{r-1}\vec e\cdot \nabla\phi$. With the help of the determination of
$H^\phi_2(T)$ in Lemma \ref{lm:koszul-phi}, this gives the existence of a homogeneous element $\vec G\in T^3$ of degree $n-2$, and homogeneous polynomials
$B, D\in \C[y,z]$ and $A\in \C[y]$ such that
$$
\vec F= \frac{\alpha}{3}\phi^{r-1}\vec e +\vec G\times \nabla\phi +
\left(\begin{smallmatrix}xD\\xA+B\\-2zD\end{smallmatrix}\right).
$$ We now compute the divergence of $\vec F$:
$$
\mathrm{Div}(\vec F) = \alpha r\phi^{r-1} +(\nabla\times \vec G)\cdot \nabla\phi -D
+xA' +\frac{\partial B}{\partial y} -2z\frac{\partial D}{\partial z}.
$$
The equation (\ref{eq:HP2}) then becomes
$$
0= (\nabla\times \vec G)\cdot \nabla\phi -D +xA' +\frac{\partial B}{\partial y}
-2z\frac{\partial D}{\partial z},
$$
which shows that $x$ divides the polynomial $\frac{\partial B}{\partial y}
-2z\frac{\partial D}{\partial z}-D\in \C[y,z]$. This implies that $
\frac{\partial B}{\partial y} =D +2z\frac{\partial D}{\partial z}
$
and $0= (\nabla\times \vec G)\cdot \nabla\phi +xA' $, so that we also have $A'\in
(xT+zT)\cap \C[y]=\lbrace 0\rbrace$. We have obtained that $A=\beta\in \C$ is a
constant and $$
(\nabla\times \vec G)\cdot \nabla\phi=0.
$$
Lemma \ref{lm:eq-phi-rot} leads to the existence of  homogeneous polynomials
$G, H\in T$,  and $C, P\in \C[y,z]$, of respective degrees equal to $n-3$ and $n-4$,
satisfying  %
 \begin{equation}\label{eq:eqCPu}
 \frac{\partial C}{\partial y} = P+2z\frac{\partial P}{\partial z}
 \end{equation}
 and $\gamma\in \C$, $k\in \mathbb N$, such that
 \begin{equation*}
 \vec G= \nabla G+H\nabla \phi
+\left(\begin{smallmatrix}z\\0\\-x\end{smallmatrix}\right)C+\left(\begin{smallmatrix}2yz\\-3xz\\xy\end{smallmatrix}\right)P
+\gamma \, \phi^k\left(\begin{smallmatrix}xz\\0\\-x^2\end{smallmatrix}\right).
 \end{equation*}
 Now, this enables us to write:
 $$
 \vec F= \frac{\alpha}{3}\phi^{r-1}\vec e+\nabla G\times \nabla\phi  +
\left(\begin{smallmatrix}xD\\B\\-2zD\end{smallmatrix}\right)
 +\left(\begin{smallmatrix}3x^3z\, P\\3x^2z \, C\\-6x^2z^2\, P\end{smallmatrix}\right)
  +(-3\gamma \, x\phi^{k+1}+\beta
x)\left(\begin{smallmatrix}0\\1\\0\end{smallmatrix}\right).
 $$
Let us now fix a homogeneous polynomial $P_1\in \C[y,z]$ verifying $\frac{\partial
P_1}{\partial y}=P$. Then the equation (\ref{eq:eqCPu}) gives $
\frac{\partial}{\partial y}\left(C-P_1-2z\frac{\partial P_1}{\partial z}\right)=0,
$
i.e., $C-P_1-2z\frac{\partial P_1}{\partial z}$ is a homogeneous polynomial in
$\C[z]$, of degree $n-3$, which means that there exists $\eta\in \C$ such that
$C=P_1+2z\frac{\partial P_1}{\partial z}+\eta z^{n-3}$. Now compute
\begin{eqnarray*}
\delta_3\left(-3xz P_1 +\frac{-3\eta}{2n-5}xz^{n-2}\right) &=& -3\nabla\left(xz P_1
+\frac{\eta}{2n-5}xz^{n-2}\right)\times \nabla\phi\\
&=& 3 \left(\begin{smallmatrix}z\, P_1+\frac{\eta}{2n-5}z^{n-2}\\xz \frac{\partial
P_1}{\partial y}\\x P_1+xz\frac{\partial P_1}{\partial z}
+\frac{\eta(n-2)}{2n-5}xz^{n-3}\end{smallmatrix}\right)\times
\left(\begin{smallmatrix}2xz\\0\\x^2\end{smallmatrix}\right)\\
&=&  3\left(\begin{smallmatrix}x^3z P\\x^2zP_1+2x^2z^2\frac{\partial P_1}{\partial
z}+\eta x^2z^{n-2}\\-2x^2z^2P\end{smallmatrix}\right)\\
&=&  \left(\begin{smallmatrix}3x^3z P\\3x^2z\, C\\-6x^2z^2P\end{smallmatrix}\right).
\end{eqnarray*}
This implies:
$\vec F \in \mathrm{Im}(\delta_3) +\C[\phi]\vec e  +
\left(\begin{smallmatrix}xD\\B\\-2zD\end{smallmatrix}\right)  +x\C[\phi]\left(\begin{smallmatrix}0\\1\\0\end{smallmatrix}\right)$. Conversely, it
is straightforward to see that an element of this space lies in the kernel of
$\delta_2$. In other words,
$$
HP_2(T) = \C[\phi]\vec e  +x\C[\phi]\left(\begin{smallmatrix}0\\1\\0\end{smallmatrix}\right)
+ \left\lbrace\left(\begin{smallmatrix}xD\\B\\-2zD\end{smallmatrix}\right)\mid D,
B\in \C[y,z]\scriptstyle \hbox{ satisfying } \frac{\partial B}{\partial y} =
D+2z\frac{\partial D}{\partial z} \right\rbrace.
$$
Let us show that the previous sum is a direct one. To do this, let us consider a
homogeneous element of this sum. Let $n\in \mathbb N$, $\alpha, \beta\in \C$, and
homogeneous polynomials $D, B\in \C[y,z]$ of respective degrees equal to $3n$ and
$3n+1$ satisfying $\frac{\partial B}{\partial y} = D+2z\frac{\partial D}{\partial
z}$ and a homogeneous $G\in T$ of degree $3n$, such that
$$
\alpha \phi^n\vec e + \beta
x\phi^n\left(\begin{smallmatrix}0\\1\\0\end{smallmatrix}\right) +
\left(\begin{smallmatrix}xD\\B\\-2zD\end{smallmatrix}\right) = \nabla
G\times\nabla\phi = \left(\begin{smallmatrix}-x^2\frac{\partial G}{\partial y}\\-2xz\frac{\partial
G}{\partial z}+x^2\frac{\partial G}{\partial x}\\2xz\frac{\partial G}{\partial
y}\end{smallmatrix}\right) .
$$
Computing the inner product of this identity with $\nabla\phi$ leads to $3\alpha
\phi^{n+1}=0$, so that $\alpha=0$.
Moreover, this gives $xD = -x^2\frac{\partial G}{\partial y}$, so that $D\in  xT\cap
\C[y,z]=\lbrace 0\rbrace$ and $\frac{\partial G}{\partial y}=0$, i.e., $G\in
\C[x,z]$.
Moreover, the second row of the previous equation implies that $B\in xT\cap
\C[y,z]=\lbrace 0\rbrace$, and $B=0$.
It remains to show that $\beta=0$, while we have $\beta
x\phi^n\left(\begin{smallmatrix}0\\1\\0\end{smallmatrix}\right) =\nabla G\times\nabla\phi $.
According to Lemma \ref{lm:phidivides}, $G$ is in $\mathbb C[\phi]$ and $\beta =0$.
\\
We then have obtained
$$
HP_2(T) = \C[\phi]\vec e \;  \oplus\;
x\C[\phi]\left(\begin{smallmatrix}0\\1\\0\end{smallmatrix}\right)
\;\oplus\;
\left\lbrace\left(\begin{smallmatrix}xD\\B\\-2zD\end{smallmatrix}\right)\mid D, B\in
\C[y,z]: \scriptstyle\frac{\partial B}{\partial y} = D+2z\frac{\partial D}{\partial
z} \right\rbrace,
$$
so that it remains to show that %
\begin{eqnarray}\label{eq:directsumHP2}
\lefteqn{\left\lbrace\left(\begin{smallmatrix}xD\\B\\-2zD\end{smallmatrix}\right)\mid
D, B\in \C[y,z]: \frac{\partial B}{\partial y} = D+2z\frac{\partial D}{\partial z}
\right\rbrace
=}\nonumber\\
&& \qquad \qquad \qquad z\mathbb
C[z]\left(\begin{smallmatrix}0\\1\\0\end{smallmatrix}\right)\oplus
\renewcommand{\arraystretch}{0.7}
  \bigoplus_{\begin{array}{c}\scriptstyle n\in \mathbb N\\\scriptstyle 0\leq k\leq
n\end{array}}\mathbb C\left(\begin{smallmatrix}(k+1) x\\ (2(n-k)+1)\,y\\
-2(k+1)\,z\end{smallmatrix}\right)y^{k}z^{n-k}.
\end{eqnarray}
To do this, let us consider $n\in \mathbb N$ and homogeneous polynomials $D, B\in
\C[y,z]$ of respective degrees $n$ and $n+1$, satisfying $\frac{\partial B}{\partial
y} = D+2z\frac{\partial D}{\partial z}$. We have already seen in the proof of
Corollary \ref{cr:equation} that this implies the existence of complex numbers
$a_k\in \C$, $0\leq k\leq n$ and $b_0\in \C$, such that $$
D= \sum_{k=0}^{n} a_k\,  y^k z^{n-k}, \hbox{ and } B =
b_0z^{n+1}+\sum_{k=1}^{n+1}{\scriptstyle\frac{2(n-k)+3}{k}}a_{k-1} \,y^k z^{n+1-k}.
$$
This gives %
\begin{eqnarray*}
\left(\begin{smallmatrix}xD\\ B\\ -2z D\end{smallmatrix}\right) &=& b_0z^{n+1}\left(\begin{smallmatrix}0\\ 1\\ 0\end{smallmatrix}\right) + \sum_{k=0}^n
\frac{a_k}{k+1} \left(\begin{smallmatrix} (k+1)\,xy^kz^{n-k}\\
(2(n-k)+1)\,y^{k+1}z^{n-k}\\ -2(k+1)\,y^kz^{n+1-k}\end{smallmatrix}\right)\\
&\in& z\mathbb C[z]\left(\begin{smallmatrix}0\\1\\0\end{smallmatrix}\right)\oplus
\renewcommand{\arraystretch}{0.7}
  \bigoplus_{k=0}^n\mathbb C\left(\begin{smallmatrix}(k+1) x\\ (2(n-k)+1)\,y\\
-2(k+1)\,z\end{smallmatrix}\right)y^{k}z^{n-k}.
\end{eqnarray*}
Notice that the previous sum is a direct one, because of degree reasons. This
allows us to conclude that (\ref{eq:directsumHP2}) holds, which finishes the
determination of $HP_2(T)$.
\Edm\\
\\
\Rm
Recall that the Poisson structure which equips the algebra $T$ is unimodular, i.e., its modular class (see
\cite{wein:mod}) vanishes (here, even the curl vector field (see \cite{dufhar:curl}) is zero), which implies that there is a duality between the Poisson cohomology and the
Poisson homology of the Poisson algebra $T$:
$
HP^{\bullet}(T)æ\simeq HP_{3-\bullet}(T).
$
\\

Now, we show in an elementary way that each Poisson cycle can be lifted to a Koszul cycle. By definition, a \emph{Koszul cycle} is a cycle of the filtered complex (\ref{shift}). %
\Bpo \label{lifting}
We keep the notations of the beginning of the section for the algebra $B$ and the Poisson algebra $(T, \lbrace\cdot, \cdot\rbrace)$.
Let us consider the filtration $F$ on the algebra $B$, given by the degree in $y$. For every Poisson cycle $X$, there exists a Koszul cycle (i.e. a cycle of the filtered complex (\ref{shift})) $\tilde X$ such that $gr_F(\tilde X) \simeq X$.
\Epo
\Bdm Proposition \ref{po:poissonhomol} gives a basis for each Poisson homology
vector space of $T$. For every Poisson boundary $\delta_k(X)$, we have seen in
Proposition \ref{po:gr} that $gr_F(\tilde d_k(Y))=\delta_k(X)$, where $Y$ is the
element $X$, viewed in the algebra $B$ and written in the basis $(x^iy^jz^k)$,
$i,j,k\in \mathbb N$. So that $\tilde X:=\tilde d_k(Y)$ is a Koszul cycle satisfying $gr_F(\tilde X) \simeq \delta_k(X)$. This implies that it suffices to show that
each element $X$ of the bases of the Poisson homology spaces given in Proposition \ref{po:poissonhomol}, can be lifted to a
Koszul cycle, i.e., for each element $X$ of the bases, we will give a Koszul
cycle $\tilde X$, satisfying $gr_F(\tilde X)=X$.
Notice that we will use here the identifications explained before: $B(k(c(w))\simeq
B$, $BR_B\simeq B^3$, $BV_B\simeq B^3$, and $\Omega^3(T)\simeq T$,
$\Omega^2(T)\simeq T^3$, $\Omega^1(T)\simeq T$ and $\Omega^0(T)\simeq T$.\\

\noindent \textit{Lifting of the Poisson $0$-cycles}.\\
As every element of $B$ is a Koszul $0$-cycle (and similarly for the Poisson
$0$-cycles), every Poisson $0$-cycle can be lifted to a Koszul $0$-cycle.\\

\noindent\textit{Lifting of the Poisson $3$-cycles}.\\
According to Proposition \ref{po:poissonhomol}, $HP_3(T)\simeq \C[\phi]=\C[x^2z]$.
For each $c\in\C$ and $k\in \mathbb N$, the definition of $\tilde d_3$ and lemma \ref{lm:form} clearly lead to
$\tilde d_3(cx^{2k}z^k\otimes c(w))= 0$, so that $cx^{2k}z^k\otimes c(w)$, which is
identified to $cx^{2k}z^k$, is a Koszul $3$-cycle satisfying
$gr_F(cx^{2k}z^k)=cx^{2k}z^k$.\\

\noindent\textit{Lifting of the Poisson $1$-cycles}.\\
 According to Proposition \ref{po:poissonhomol}, the space $HP_1(T)$ is generated as
a $\C$-vector space by the following elements: $A_{k}:=
(x^2z)^k\left(\begin{smallmatrix}xz\\0\\ -x^2\end{smallmatrix}\right)$, $B_{r} :=
z^r\left(\begin{smallmatrix}z^2\\0\\-xz\end{smallmatrix}\right)$, $\vec u_{n,k} =
\left(\begin{smallmatrix}(2n+3)\,y^kz^{n+2-k}\\-3k\,xy^{k-1}z^{n+2-k}\\(-2n+3(k-1))\,xy^kz^{n+1-k}\end{smallmatrix}\right)$,
$\vec v_{m,
s}:=\left(\begin{smallmatrix}0\\sy^{s-1}z^{m-s}\\(m-s)y^sz^{m-1-s}\end{smallmatrix}\right)$,
and $\vec
w_{p}:=\left(\begin{smallmatrix}y^p\\p\,xy^{p-1}\\0\end{smallmatrix}\right)$,
 where $n, p\in \mathbb N$, $m\in \mathbb N^*$ and $ 1\leq k\leq n+1$, $0\leq s\leq m$.

Now, by definition of $\tilde d_1$ and because $xz=zx$ in $B$, it is clear that:
$$
\tilde d_1(x^{2k+1}z^{k+1} \otimes x - x^{2k+2}z^k\otimes z)= 0, \quad \hbox{ and }
\quad \tilde d_1(z^{r+2}\otimes x-xz^{r+1}\otimes z) = 0,
$$
so that $\tilde A_k := x^{2k+1}z^{k+1} \otimes x - x^{2k+2}z^k\otimes z$ and $\tilde
B_r  := z^{r+2}\otimes x-xz^{r+1}\otimes z$ are Koszul $1$-cycles satisfying
$gr_F(\tilde A_k)=A_k$ and $gr_F(\tilde B_r)=B_r$.

Next, we let
$$
\begin{array}{l}
\tilde U_{n,k} :=\\
\\
 (2n+3)\,y^kz^{n+2-k}\otimes x-3k\,xy^{k-1}z^{n+2-k}\otimes y
+{(-2n+3(k-1))}\,xy^kz^{n+1-k}\otimes z\\
\\
\left\lbrace\begin{array}{l}
 -\sum\limits_{\ell=0}^{k-2} \frac{k!}{\ell!}
\;\frac{a_{n,k,\ell}}{2(n+2-k)-(k-\ell)}\; \mathcal{X}_{k-\ell, \ell, n+2-k},\quad
\hbox{ if } \; 3k-2(n+2)<0;\\
 \\
 \\
 + \sum\limits_{\ell=0}^{\ell_{(n,k)}-1} \frac{k!}{\ell!}
\frac{\left(a_{n,k,\ell_{(n,k)}}-a_{n,k,\ell}\right)}{\ell-\ell_{(n,k)}}
\mathcal{X}_{k-\ell, \ell,
n+2-k}-\sum\limits_{\ell=\ell_{(n,k)}+1}^{k-2}\frac{k!}{\ell !} \frac{a_{n,k,\ell}}{\ell-\ell_{(n,k)}}\mathcal{X}_{k-\ell, \ell, n+2-k}\\  \\
   -\frac{k!}{(\ell_{(n,k)}+1)!} a_{n,k,\ell_{(n,k)}}\;x^{2n+3-2k}
y^{3k-3-2n}z^{n+2-k}\otimes x, \quad \hbox{ if } \; 3k-2(n+2)\geq 0,\\
 \end{array}\right.
\end{array}
$$
where for $1\leq k\leq n+1$ and $0\leq \ell\leq k-2$, $a_{n,k,\ell}:=-2n-6+6k-3\ell+2n(k-\ell)-3k(k-\ell)$, for all $a,b,c\in \mathbb N$, $\mathcal{X}_{a,b,c}:=x^ay^b z^c\otimes y - 2c\;x^{a}y^b z^c\otimes x$ and  if $3k-2(n+2)\geq 0$, $\ell_{(n,k)}:=3k-4-2n$.
Then, we have $gr_F(\tilde U_{n,k})=\vec u_{n,k}$ and lemma \ref{lm:form} allows one to verify that $\tilde U_{n,k}$ is a
Koszul $1$-cycle.

Now, let $m\in \mathbb N^*$, $0\leq s\leq m$ and consider
$$
\begin{array}{l}
\tilde V_{m,s} := s\, y^{s-1}z^{m-s}\otimes y+(m-s)\, y^sz^{m-1-s}\otimes z \\
\\
\qquad \qquad \qquad \qquad+
\sum\limits_{\ell=0}^{s-2}\;\frac{(m-s)}{2(m-s)-1}\;\frac{s!}{\ell ! }
\;x^{s-\ell-1}y^\ell z^{m-s}\otimes y.
\end{array}
$$
Then we have $gr_F(\tilde V_{m,s})=\vec v_{m,s}$ and it is straightforward, using
lemma \ref{lm:form}, to verify that $\tilde V_{m,s}$ is
a Koszul $1$-cycle.

Finally, let $p\in \mathbb N$ and consider $$
\tilde W_p := \sum_{k=0}^p \frac{p!}{k!}\, x^{p-k} y^k\otimes x + \sum_{k=0}^{p-1}
\frac{p!}{k!} x^{p-k} y^k\otimes y.
$$
It is clear that $gr_F(\tilde W_p)=\vec w_p$ and moreover,
using lemma \ref{lm:form}, we obtain that $\tilde
d_1(\tilde W_p)=0$, i.e., $\tilde W_p$ is a Koszul $1$-cycle.\\

\noindent\textit{Lifting of the Poisson $2$-cycles}.\\
 According to Proposition \ref{po:poissonhomol}, the space $HP_2(T)$ is generated as
a $\C$-vector space by the following elements: $C_r := (x^2z)^r \left(\begin{smallmatrix}x\\y\\z\end{smallmatrix}\right)$, $D_s :=
x(x^2z)^s\left(\begin{smallmatrix}0\\1\\0\end{smallmatrix}\right)$, $E_t
:=z^{t+1}\left(\begin{smallmatrix}0\\1\\0\end{smallmatrix}\right)$, $\vec o_{n,k} :=
\left(\begin{smallmatrix}(k+1) xy^{k}z^{n-k}\\ (2(n-k)+1)\,y^{k+1}z^{n-k}\\
-2(k+1)\,y^{k}z^{n-k+1}\end{smallmatrix}\right)$,
where $r, s, t, n\in \mathbb N$ and $0\leq k\leq n$.

First, if $r\in \mathbb N$, then, using the definition of $\tilde d_2$ and lemma \ref{lm:form}, we obtain easily:
$$
\tilde d_2(x^{2r+1}z^r \otimes r_1 + x^{2r}yz^r \otimes r_2 +x^{2r}z^{r+1} \otimes
r_3) = 0,
$$
so that $\tilde C_r := x^{2r+1}z^r \otimes r_1 + x^{2r}yz^r \otimes r_2+x^{2r}z^{r+1} \otimes r_3$ is a Koszul $2$-cycle satisfying $gr_F(\tilde C_r)=C_r$.

Moreover, if $s, t\in \mathbb N$, then, because $xz=zx$ in $B$ and by definition of
$\tilde d_2$, we have
$$
\tilde d_2(x^{2s+1}z^s \otimes r_2)=0, \qquad \hbox{ and }\qquad \tilde d_2(z^{t+1}\otimes r_2)=0, $$
so that $\tilde D_s:= x^{2s+1}z^s \otimes r_2$ and $\tilde E_t := z^{t+1}\otimes
r_2$ are Koszul $2$-cycles and satisfy $gr_F(\tilde D_s)=D_s$ and $gr_F(\tilde
E_t)= E_t$.

Now, let $ n\in \mathbb N$ and $0\leq k\leq n$ and %
\begin{eqnarray*}
\lefteqn{\tilde O_{n,k} := (k+1)\, xy^{k}z^{n-k}\otimes r_1+
(2(n-k)+1)\,y^{k+1}z^{n-k}\otimes r_2}\\
&&\qquad \qquad+\left( -2(k+1)\,y^{k}z^{n-k+1}-\sum_{j=0}^{k-1}\frac{(k+1)!}{j!}\,
x^{k-j} y^jz^{n+1-k}\right)\otimes r_3.
\end{eqnarray*}
Then, using once more lemma \ref{lm:form}, it is straightforward to verify that $\tilde d_2(\tilde
O_{n,k})=0$. Moreover, we have of course, $gr_F(\tilde O_{n,k})=\vec o_{n,k}$, which
finishes the proof.
\Edm\\
\\
\emph{Proof of Theorem \ref{calc}}. Actually, it remains to prove that the Hochschild
homology of the algebra $B$ is isomorphic to the Poisson homology of $T$ obtained in Proposition \ref{po:poissonhomol}. Following the same method as in~\cite{vdb:Kth}, we use the Brylinski spectral sequence of the almost commutative algebra $B$~\cite{bry:poi,kas:homcy}. Denote by $C$ the filtered complex (\ref{shift}) and denote by $(F_pC)_{p\in \mathbb{Z}}$ its filtration. The complex $F_pC$ is the following %
\begin{equation}\label{FpC}
0 \longrightarrow F_{p+3}(B(\mathbb{C}c(w)))\stackrel{\tilde d_3}{\longrightarrow}  F_{p+2}(BR_B)\stackrel{\tilde d_2}{\longrightarrow} F_{p+1}(BV_B) \stackrel{\tilde d_1}{\longrightarrow} F_p(B) \longrightarrow 0
\end{equation}
where $F_p$ denotes the filtration of the total degree. Let us consider the spectral sequence associated to the filtered complex $C$ (Section 5.4 in~\cite{weib:homo}). The term $E^0$ of this spectral sequence is the graded complex naturally associated to the filtered complex $C$. By Proposition \ref{po:gr}, $E^0$ is isomorphic to the Poisson complex of the Poisson algebra $T$ (actually, the proof of Proposition \ref{po:gr} provides an explicit isomorphism). Since the filtration $F$ of the complex $C$ is increasing, exhaustive and bounded below ($F_{-1}B = 0$), the spectral sequence converges to $H_{\bullet}(C)$ (Theorem 5.5.1.2 in~\cite{weib:homo}):
\begin{equation}\label{specseq}
E^1_{pq}=H_{p+q}(F_pC/F_{p-1}C) \Longrightarrow H_{p+q}(C).
\end{equation}
Thus, in order to conclude that the Hochschild
homology of the algebra $B$ is isomorphic to the Poisson homology of $T$, it is sufficient to prove the following.
\Bpo
The spectral sequence associated to the filtered complex $C$ degenerates at $E^1$.
\Epo
\Bdm
We apply a standard criterion for degeneration of spectral sequences (Lemma 5.2 in~\cite{vdb:Kth}) with $r=1$ (we use the notation of~\cite{vdb:Kth}). This criterion consists in proving that the natural edge map
$$
\phi^1_p : H_{\bullet}(F_pC) \longrightarrow E^1_p
$$
is surjective for any $p$. Since the term $E^1$ is isomorphic to the Poisson homology of $T$, surjectivity is given by Proposition \ref{lifting}.  \Edm\\

Actually the proof of Proposition \ref{lifting} provides an explicit section of the edge map $\phi^1_p$. Consequently we have an explicit isomorphism from $HP_{\bullet}(T)$ to $HH_{\bullet}(B)$ described as follows: Proposition \ref{po:poissonhomol} gives an explicit family $(c)$ of Poisson cycles $c$ such that the family $([c])$ of their classes form a basis of the space $HP_{\bullet}(T)$, and the proof of Proposition \ref{lifting} lifts each Poisson cycle $c$ to an explicit Koszul cycle $\tilde c$, so that the Hochschild classes $[\tilde c]$ form a basis of $HH_{\bullet}(B)$.

Since $B$ is 3-Calabi-Yau (Theorem \ref{calyau}), we deduce Hochschild cohomology of $B$ from Theorem \ref{calc}: $HH^{\bullet}(B) \cong HH_{3-\bullet}(B)$.


\begin{thebibliography}{99}


\bibitem{as:regular} M. Artin, W. F. Schelter, Graded algebras of global dimension
3, \emph{Adv.
Math.} \textbf{66} (1987), 171-216.
\bibitem{bafr:kalg} J. Backelin, R. Fr\"{o}berg, Koszul algebras, Veronese subrings
and rings with
linear resolutions, \emph{Rev. Roum. Math. Pures Appli.} \textbf{30} (1985), 85-97.
\bibitem{rb:gera} R. Berger, Gerasimov's theorem and $N$-Koszul algebras, \emph{J.
London Math. Soc.} \textbf{79} (2009), 631-648.
\bibitem{bdvw:homog} R. Berger, M. Dubois-Violette, M. Wambst, Homogeneous algebras,
\emph{J. Algebra} \textbf{261} (2003), 172-185.
\bibitem{rbrt:pbw} R. Berger, R. Taillefer, Poincar\'e-Birkhoff-Witt deformations of
Calabi-Yau algebras, \emph{J. Noncommut. Geom.} \textbf{1} (2007), 241-270.
\bibitem{bock:graded} R. Bocklandt, Graded Calabi-Yau algebras of dimension 3,
\emph{J. Pure Appl. Algebra} \textbf{212} (2008), 14-32. \bibitem{bg:aqg} K. A. Brown, K. R. Goodearl, \emph{Lectures on algebraic quantum
groups}, Advanced Courses in Mathematics CRM Barcelona, Birkha\"user, 2002. \bibitem{bry:poi} J.L. Brylinski, A differential complex for Poisson manifolds,
\emph{J. Differential Geom.} \textbf{28} (1988), 93-114.
\bibitem{cktb:def} A. Cattaneo, B. Keller, C. Torossian, A. Brugui\`eres,\emph{D\'eformation, quantification, th\'eorie de Lie},
Panoramas et Synth\`eses, SMF, 2005.
\bibitem{cw:bili} B. Corbas, G. D. Williams, \emph{Bilinear forms over an algebraically closed field}, J. Pure Appl. Algebra \textbf{165} (2001), 255-266.
\bibitem{doik:congr} D. Z. Dokovic, K. D. Ikramov, A square matrix is congruent to its transpose, \emph{J. Algebra} \textbf{257} (2002), 97-105.
\bibitem{dv:multi} M. Dubois-Violette, Multilinear forms and graded algebras, \emph{J. Algebra} \textbf{317} (2007), 198-225. \bibitem{dufhar:curl} J.P. Dufour, A. Haraki, Rotationnels et structures de Poisson quadratiques, \emph{C. R. Acad. Sci. Paris S\'er. I Math.} \textbf{312} (1991), 137--140.
\bibitem{vg:cy} V. Ginzburg, Calabi-Yau algebras, arXiv:math.AG/0612139.
\bibitem{hot:orbi} G. Halbout, J.-M. Oudom, X. Tang, Deformations of orbifolds with
noncommutative linear Poisson structures, \emph{Int. Math. Res. Not.} (2011), 1-39. \bibitem{kas:homcy} C. Kassel, L'homologie cyclique des alg\`ebres enveloppantes, \emph{Invent. Math.} \textbf{91} (1988), 221-251.
\bibitem{bk:compl} B. Keller, Deformed Calabi-Yau completions (With an appendix by Michel Van den Bergh), \emph{Journal f\"ur die Reine und Angewandte Mathematik}, \textbf{654} (2011), 125-180. \bibitem{ks:qg} A. Klimyk, K. Schm\"udgen, \emph{Quantum groups and their representations},
Texts and Monographs in Physics, Springer, 1997.
\bibitem{kont:quant} M. Kontsevich, Deformation quantization of Poisson manifolds, \emph{Lett. Math. Phys.} \textbf{66} (2003), 157-216.
\bibitem{lbsvdb:central} L. Le Bruyn, S. P. Smith, M. Van den Bergh, Central
extensions of three-dimensional Artin-Schelter regular algebras, \emph{Math. Z.} \textbf{222} (1996), 171-212.
\bibitem{man:kalg} Y. I. Manin, \emph{Quantum groups and non-commutative geometry},
CRM, Universit\'e de Montr\'eal, 1988.
\bibitem{nic:cubic} N. Marconnet, Homologies of cubic Artin-Schelter regular
algebras, \emph{J. Algebra} \textbf{278} (2004), 638-665. \bibitem{mon:quad} P. Monnier, Formal Poisson cohomology of quadratic Poisson
structures, \emph{Lett. Math. Phys.} \textbf{59} (2002), 253-267. \bibitem{pel:four} S. Pelap, Poisson (co)homology of polynomial Poisson algebras in
dimension four: Sklyanin's case, \emph{J. Algebra} \textbf{322} (2009), 1151-1169.
\bibitem{pic:hom} A. Pichereau, Poisson (co)homology and isolated singularities, \emph{J. Algebra} \textbf{299} (2006), 747-777.
\bibitem{popo:quad} A. Polishchuk, L. Positselski, \emph{Quadratic algebras},
University Lecture Series 37, AMS, 2005.
\bibitem{ser:classif} V. V. Sergeichuk, Classification problems for systems of forms and linear mappings, \emph{Math. USSR Izvestiya} \textbf{31} (1988), 481-501.
\bibitem{smith:octo} S. P. Smith, A 3-Calabi-Yau algebra with $G_2$ symmetry
constructed from the octonions, arXiv:1104.3824v1.
\bibitem{sz:gro} D. R. Stephenson, J. J. Zhang, Growth of graded noetherian rings,
\emph{Proc. Amer. Math. Soc.} \textbf{125} (1997), 1593-1605.
\bibitem{msa:steiner} M. Su\'{a}rez-Alvarez, 3-Calabi-Yau algebras from Steiner
triple systems,
preprint May 2011.
\bibitem{vdb:Kth} M. Van den Bergh, Non-commutative homology of some
three-dimensional quantum spaces, \emph{K-Theory} \textbf{8} (1994), 213-220.
\bibitem{vdb:ober} M. Van den Bergh, Introduction to superpotentials, \emph{Oberwolfach Rep.} \textbf{2} (2005), 396-398. \bibitem{wam:kosz} M. Wambst, Complexes de Koszul quantiques, \emph{Ann. Inst.
Fourier} \textbf{43} (1993), 1089-1156.  \bibitem{weib:homo} C. A. Weibel, \emph{An introduction to homological algebra},
Cambridge University
Press, 1994.
\bibitem{wein:mod} A. Weinstein, The modular automorphism group of a Poisson manifold, \emph{Journal of Geometry and Physics} \textbf{23} (1997), 379--394. \bibitem{jz:nnoeth} J. J. Zhang, Non-noetherian regular rings of dimension 2,
\emph{Proc. Amer. Math. Soc.} \textbf{126} (1998), 1645-1653.
\end{thebibliography}
\end{document}